\documentclass[a4paper, 10pt]{article}
\usepackage{amssymb, amsthm, amsmath, latexsym}
\usepackage{amsfonts}

\usepackage{multicol}
\usepackage{euscript}
\usepackage{amssymb}
\usepackage{hyperref}
\usepackage{amsmath}
\usepackage{amsthm}
\usepackage[dvips]{graphicx}
\usepackage{psfrag}
\newcommand{\unpo}{\mathcal{O}(1)}
\newcommand{\ee}{\varepsilon}

\newcommand{\zz}{\zeta}

\newcommand{\uu}{\Upsilon}

\newtheorem{teo}{Theorem}[section]
\newtheorem{pro}{Proposition}[section]
\newtheorem{lem}{Lemma}[section]

\theoremstyle{definition}
\newtheorem{rem}{Remark}[section]

\newtheorem{hyp}{Hypothesis}

\newtheorem{say}{Definition}[section]

\begin{document}
\numberwithin{equation}{section}
\bibliographystyle{plain}

\title{Invariant manifolds for a singular ordinary differential equation}
\author{ \it Stefano Bianchini \\ \footnotesize SISSA, via Beirut 2-4, 34014, Trieste, Italy \\ \footnotesize email: bianchin@sissa.it 
 \\  \small and \\  \it Laura V. Spinolo   \\   \footnotesize
Centro di Ricerca Matematica Ennio De Giorgi, Scuola Normale Superiore,  
Piazza dei Cavalieri 3, 
56126 Pisa, 
Italy \\ \footnotesize email: laura.spinolo@sns.it}     
         
\date{}
\maketitle
\begin{abstract}
           We study the singular ordinary differential equation
           \begin{equation}
           \label{e:i:abstract}
                      \frac{d U}{d t} = \frac{1}{\zeta (U)} \phi_s (U) + \phi_{ns} (U),
           \end{equation}
           where $U \in \mathbb R^N$, the functions $\phi_s \in \mathbb R^N $ and $\phi_{ns} \in \mathbb R^N $ are of class $\mathcal C^2$ and 
           $\zeta $ is  a real valued $\mathcal C^2$ function. The equation is singular in the sense that $\zeta (U)$ can attain the value $0$. 
           We focus on the solutions of \eqref{e:i:abstract} that belong to a small neighbourhood of a point $\bar U$ such that 
           $\phi_s (\bar U) = \phi_{ns} (\bar U) = \vec 0$, $\zz (\bar U) =0$. We investigate the existence of manifolds that are locally invariant for \eqref{e:i:abstract} and that contain orbits 
           with a suitable prescribed asymptotic behaviour. Under suitable hypotheses on the set $\{ U: \; \zz (U) = 0 \}$, we extend to the case of the singular ODE \eqref{e:i:abstract} the definitions of  center manifold,  center stable manifold and of uniformly stable manifold. We prove that the solutions of \eqref{e:i:abstract} lying on each of these manifolds are regular: this is not trivial since we provide examples showing that, in general, a solution of \eqref{e:i:abstract} is not continuously differentiable.   
           Finally, we show a decomposition result for a center stable and for the uniformly stable manifold. 
           
           An application of our analysis concerns the study of the viscous profiles with small total variation for a class of mixed hyperbolic-parabolic systems in one space variable. Such a class includes the compressible Navier Stokes equation.
\vspace{0.3cm}

\noindent {\bf Key words:} singular ordinary differential equation, stable manifold, center manifold, invariant manifold.  
 \end{abstract}
\section{Introduction}
\label{s:intro}
In this work we study the singular ordinary differential equation
\begin{equation}
\label{e:intro:theq}
             \frac{d U}{d t} = \frac{1}{\zeta (U)} \phi_s (U) + \phi_{ns} (U).
\end{equation}
In the previous expression, $U \in \mathbb R^N$ and  the functions $\phi_s$ and $\phi_{ns}$ are $\mathcal C^2$ (continuously differentiable with continuously differentiable derivatives) and take values into $\mathbb R^N$. The function $\zz $ is as well regular and it takes real values. We say that the equation is singular because $\zz (U)$ can attain the value $0$.

Equation \eqref{e:intro:theq} is related to a class of problems studied in singular perturbation theory. Consider system
\begin{equation}
\label{e:i:spt}
   \left\{
   \begin{array}{ll}
                \ee d \mathbf{x} / dt = f (\mathbf{x}, \, \mathbf{y}, \, \ee) \\
                d \mathbf{y} / dt = g (\mathbf{x}, \, \mathbf{y}, \, \ee), \\
    \end{array}            
   \right.
\end{equation}
where $\mathbf{x}$ and $\mathbf{y}$ are vector valued functions, $\ee$ is a parameter. In singular perturbation theory one is typically concerned with the limit $\ee \to 0$ and with the corresponding behaviour of the solution $(\mathbf{x}, \, \mathbf{y})$.  Note that \eqref{e:intro:theq} can be viewed as as an extension of \eqref{e:i:spt}, in the sense that \eqref{e:i:spt} can be written in the form \eqref{e:intro:theq}: in this case, the singularity $\zz (U)$ in \eqref{e:intro:theq} is identically equal to $\ee$ and hence $d \zz / dt =0$. 

Being the literature concerning \eqref{e:i:spt} extremely wide, it would be difficult to give an overview here. Consequently, we just refer to the notes by Jones \cite{Jones} and to the rich bibliography contained therein. In particular, \cite{Jones} provides a nice overview of Fenichel's papers \cite{Fenichel, Fenichel_per, Fenichel_asy}. The works \cite{Fenichel, Fenichel_per, Fenichel_asy} provide several ideas and techniques exploited in the present paper. 

The main novelty of the present work is that we consider the case $\zz$ is  a nontrivial function of the unknown $U$. In particular, this means that $d \zz / dt \neq 0$ in general and hence that we have to face the possibility that $\zz\big( U (0) \big) \neq 0$,  but $\zz \big( U(t ) \big) =0$ for a finite value of $t$. This is exactly what happens in the examples \eqref{e:ex:fast} and \eqref{e:ex:slow} discussed in Section \ref{s:hyp} here. Other examples are provided in a previous work by the same authors \cite{BiaSpi:rie}, Section 2. Note that, in all these cases, there is a loss of regularity at the time $t_0$ at which $\zz \big((U (t) \big)$ reaches the value $0$, $t_0 = \min \big\{ t \in [0, \, + \infty [ : \; \zz \big(U (t ) \big) =0 \big\}$. More precisely, the first derivative $d U / dt$ either has a discontinuity or blows up at $t= t_0$.

Our goal here is to study the solutions of \eqref{e:intro:theq} that lie in a neighborhood of a point $\bar U$ such that $\phi_s (\bar U) = \phi_{ns}(\bar U) = \vec 0$ and $\zz (\bar U) = 0$.
We are concerned with the existence of invariant manifolds. More precisely, the problem is the following. 

Consider first the case of the non singular ODE
\begin{equation}
\label{e:i:nons}
          \frac{dU}{dt} = f (U)
\end{equation}
and assume that the point $\bar U$ is an equilibrium, namely $f (\bar U) = \vec 0$. In a neighbourhood of $\bar U$ one can define a \emph{center} and a \emph{center stable manifold}, which are both locally invariant for \eqref{e:i:nons}. We recall here that, loosely speaking, a center stable manifold contains the orbits of \eqref{e:i:nons} that as $ x \to + \infty$ either do not blow up or blow up more slowly than $e^{\eta t}$, where $\eta$ is a small enough constant depending on the system. More precisely, the orbits that lie on a center stable manifold are those having the asymptotic behaviour described before and solving a suitable system which, in a small neighbourhood of $\bar U$, coincides with the original one \eqref{e:i:nons}. 

Also, assume that there exists a manifold $E$ entirely constituted by equilibria of \eqref{e:i:nons} and containing $\bar U$. We are interested in the \emph{uniformly stable} manifold relative to $E$. 
By uniformly stable manifold we mean the slaving manifold that contains all the orbits that when $t \to + \infty$ decay with exponential speed to some point in $E$. Note that the uniformly stable manifold does not coincide, in general, with the classical \emph{stable} manifold. Indeed, the stable manifold contains the orbits that decay exponentially fast to the given equilibrium $\bar U$, while on the uniformly stable manifold we only require that the limit belongs to $E$. 

The existence of a center stable and of  the uniformly stable manifold can be obtained as consequence of the Hadamard-Perron Theorem, which is discussed for example in the book by Katok and Hasselblatt \cite{KHass}. 

In the present paper we prove that, under suitable hypotheses, one can extend the definition of center, center stable and of uniformly stable manifold from the case of the non singular ODE \eqref{e:i:nons} to the singular case \eqref{e:intro:theq}. The manifolds we define are all locally invariant for \eqref{e:intro:theq} and satisfy the following property: 
\begin{center}
    {\bf (P)} If $U$ is an orbit lying on the manifold and $\zz \big( U (0) \big) \neq 0$, then $\zz \big( U(t) \big) \neq 0$ for every $t$. 
\end{center}
This, in particular, rules out the losses of regularity (blow up or discontinuity in the first derivative) mentioned before. 

We proceed as follows. First, we consider the non singular ODE
\begin{equation}
\label{e:intro:theq11}
             \frac{d U}{d \tau} = \phi_s (U) + \zeta (U) \phi_{ns} (U).
\end{equation}
Formally, \eqref{e:intro:theq11} is obtained from \eqref{e:intro:theq} \emph{via} the change of variable $\tau = \tau (t)$, defined as the solution of the Cauchy problem 
\begin{equation}
\label{e:i:changev}
\left\{
\begin{array}{lll}
          \displaystyle{   \frac{d \tau}{dt} = \frac{1}{\zeta [U(t) ]} } \\
         \\
          \tau ( 0) =0.\\
\end{array}
\right.
\end{equation}
However, the function $\tau (t)$ is well defined only if $\zeta [U(t)] \neq 0 $ for every $t$. In the work we always refer to the formulation \eqref{e:intro:theq11} and we prove the existence of locally invariant manifolds satisfying property {\bf (P)}. We then show that \emph{a posteriori} the change of variable \eqref{e:i:changev} is well defined and system \eqref{e:intro:theq11} is equivalent to \eqref{e:intro:theq} on these manifolds. 

We assume that 
\begin{enumerate}
\item For every $\mathcal M^c$ center manifold of \eqref{e:intro:theq11}, the intersection between the set $\{U: \, \zeta(U) = 0 \}$ and $\mathcal M^c$ contains only equilibria.
\end{enumerate}
We then define a \emph{manifold of slow dynamics} as a center manifold of \eqref{e:intro:theq11} (any center manifold works). To simplify the exposition, in the following we \emph{fix} a manifold of slow dynamics. 
To define the \emph{manifold of fast dynamics} we assume 
\begin{enumerate} \setcounter{enumi}{1}
\item there exists a one-dimensional manifold which is transversal to the set $\{ U: \, \zz(U) = 0 \} $ and is entirely constituted by equilibria of \eqref{e:intro:theq11}. In the following, we denote by $E$ this manifold. 
\end{enumerate}
The \emph{manifold of fast dynamics} is then defined as the uniformly stable manifold of \eqref{e:intro:theq} relative to the manifold of equilibria $E$. Namely, all the fast dynamics converge exponentially fast to some equilibrium in $E$. 

The manifolds of slow and fast dynamics can be regarded as extensions to the general case of the notions of \emph{slow} and \emph{fast time scale} discussed for example in Fenichel \cite{Fenichel} in relation to system \eqref{e:i:spt}, namely to the case $\zz$ is a parameter.

We also assume that
\begin{enumerate}\setcounter{enumi}{2}
\item The singular set $\{\zeta(U)=0\}$ is invariant for both \eqref{e:intro:theq11} and for the solutions of \eqref{e:intro:theq} that lie on the manifold of the slow dynamics.
\end{enumerate}

It is not hard to show that, as a consequence of the previous assumptions, if we restrict system \eqref{e:intro:theq11} on the manifold of slow dynamics, then \eqref{e:intro:theq11} is equivalent to
\eqref{e:intro:theq}. We can thus go back to the original variable $t$ and get that the solutions of \eqref{e:intro:theq} 
lying on the manifold of the slow dynamics satisfy a non singular ODE. 
We then define a center manifold of \eqref{e:intro:theq} as a center manifold of  the system reduced on the manifold of slow dynamics (Theorem \ref{t:center}). In this way property {\bf (P)} is automatically satisfied on any center manifold  and the losses of regularity are ruled out. As before, by \emph{loss of regularity} we mean the blows up or discontinuities in the first derivative that may be exhibited by the solutions  of \eqref{e:intro:theq}, as shown by the examples  \eqref{e:ex:fast} and \eqref{e:ex:slow} in Section \ref{s:hyp} here.

To extend  to the case of the singular ODE \eqref{e:intro:theq} the definition of center stable and uniformly stable manifold we need some more work.
As mentioned before, thanks to assumption 2, there exists a manifold of equilibria transversal to the singular surface: we denoted this manifold by $E$. To define the uniformly stable manifold of \eqref{e:intro:theq} relative to $E$ we need to study the solutions of   
\eqref{e:intro:theq} which converges to an equilibrium on $E$ with exponential speed. Note that this speed can be either bounded or unbounded as $\zz \to 0$, so we are looking for a composition of both fast and slow dynamics. Also, to define a center stable manifold loosely speaking we have to study orbits that are local solutions of \eqref{e:intro:theq} and that do not blow too fast when $t \to + \infty$. Therefore, we have to deal again with a composition of slow and fast dynamics. 

In both cases (uniformly stable and center stable manifold) the analysis can be seen as an extension of the exponential splitting methods for non singular ODEs like \eqref{e:i:nons}. However, as mentioned before what \emph{a priori} can go wrong is that in the change of time scale defined by the Cauchy problem
\begin{equation}
\label{e:rescaae}
           \left\{
\begin{array}{lll}
          \displaystyle{   \frac{d \tau}{dt} = \frac{1}{\zeta [U(t) ]} } \\
         \\
          \tau ( 0) =0.\\
\end{array}
\right.
\end{equation}
some regularity is missing. 

The main result of this paper is the following (a more precise statement is given in Theorem \ref{t:app:man}):

\begin{teo} 
\label{t:intro:invman}
There is a sufficiently small constant $\delta >0$ such that the following holds. In the ball of center $\bar U$ and of radius $\delta$ in $\mathbb R^N$ 
                     one can define two  continuously differentiable manifolds $\mathcal M^s$ and $\mathcal M^{cs}$ which are both locally invariant for \eqref{e:intro:theq}. The first one, $  \mathcal M^s$,  is the uniformly stable manifold of \eqref{e:intro:theq} relative to $E$, while $\mathcal M^{cs}$ is a center stable manifold for  \eqref{e:intro:theq}. If $U$ is a solution satisfying $\zz \big( U(0) \big) \neq 0$ and lying on either $\mathcal M^s$ or $\mathcal M^{cs}$, the Cauchy problem \eqref{e:rescaae} defines a diffeomorphism $\tau: [0, \, +\infty [ \to    [0, \, +\infty [$. In other words, if we restrict to either $\mathcal M^s$ or $\mathcal M^{cs}$, then the formulations \eqref{e:intro:theq} and \eqref{e:intro:theq11} are equivalent, provided that $\zz \big( U(0) \big) \neq 0$.  In particular, property {\bf (P)}  is satisfied on both $\mathcal M^s$ and $\mathcal M^{cs}$. 
                     
                     Also, if $U(\tau)$ is a solution lying on either $\mathcal M^s$ or $\mathcal M^{cs}$, then it can be decomposed as                 
                      \begin{equation}
                      \label{e:intro:dec}
                               U( \tau) =  U_f (\tau) + U_{sl} (\tau) + U_p (\tau),
                     \end{equation}
                     where $U_{sl}$ lies on the manifold of slow dynamics and $U_f (\tau)$ is exponentially decreasing to $\vec 0$. The perturbation term $U_p (\tau)$ is small in the sense that 
                     $$
                        |U^p (\tau) | \leq k_p   \big| \zz \big( U^{\sl} (0) \big)\big| \,  \big|U_f (0) \big| e^{- c \tau / 4}
                     $$
                     for suitable positive constants $c, \; k_p >0$. 
                     \end{teo}

From the technical point of view, the key points in the proof of  Theorem \eqref{t:intro:invman} are the following two. First, we introduce a change of variables which allows us to write system
\begin{equation}
\label{e:intro:tau}
        \frac{d U}{d \tau} = \phi_s (U) + \zeta (U) \phi_{ns} (U).
\end{equation}
in a more convenient form. The precise statement is given in Proposition \ref{p:change}. This change of variables exploits many of the ideas which in the case of equation \eqref{e:i:spt} lead to the introduction of the so called \emph{Fenichel Normal Form} (see Jones \cite{Jones} and the references therein) .  Here, however, we have to rely on the assumptions 1, 2 and 3 discussed above.

The second main point in the proof of Theorem \ref{t:intro:invman} is the analysis of a family of slaving manifolds for system \eqref{e:intro:tau}.  
This analysis exploits the presence of a splitting based on exponential decay estimates and it is in the spirit of Hadamard Perron Theorem (see for example the book by Katok and Hasselblatt \cite{KHass}). Here the main results are Theorem \ref{p:unif} and Proposition \ref{p:exp:d}. Loosely speaking,  Proposition \ref{p:exp:d} tells us the following. Fix a manifold $\mathcal S$, locally invariant for \eqref{e:intro:tau} and entirely made by slow dynamics. Then there exists a slaving manifold containing orbits that decay to an orbit in $\mathcal S$ exponentially fast, with respect to the $\tau$ variable. Also, Proposition \ref{p:exp:d} ensures that any solution $U$ lying on the slaving manifold admits a decomposition like \eqref{e:intro:dec}, namely 
$$
       U( \tau) =  U_f (\tau) + U_{sl} (\tau) + U_p (\tau)
$$
where $U_{sl}$ lies on $\mathcal S$ and $U_f (\tau)$ is exponentially decreasing to $\vec 0$. The perturbation term $U_p (\tau)$ is small and disappears when $ \zz (U) = 0 $, so on the singular surface $\{ U : \; \zz (U) = 0 \}$ there is no interaction, but a complete decoupling. Actually, in the statement of Proposition \ref{p:exp:d} we consider slightly more general conditions ensuring that the interaction term disappears. However, the case $ \zz (U) = 0 $ is the one we exploit in the following. From the technical point of view, the most complicated point in the analysis is proof of the $\mathcal C^1$ regularity of the slaving manifold, since it involves studying the Frech\'{e}t  differentiability of suitable maps between Banach spaces. 
\vspace{1cm}

An application of our analysis concerns the study of the viscous profiles with small total variation for a class of mixed hyperbolic-parabolic systems in one space dimension. The connection between these viscous profiles and the singular ordinary differential equation \eqref{e:intro:theq} is discussed in \cite{BiaSpi:pro}, where we also explain what we mean by viscous profiles and by mixed hyperbolic-parabolic systems in this context. In  \cite{BiaSpi:pro} we also discuss a remark due to Fr\'{e}deric Rousset \cite{Rousset:ns} about the Lagrangian and the Eulerian formulation of the Navier Stokes equation. Loosely speaking, the connection between viscous profiles and singular ODEs like \eqref{e:intro:theq} is that the equation satisfied by the viscous profiles may be singular when the system does not satisfy a condition of block linear degeneracy defined in \cite{BiaSpi:rie}. In particular, this happens in the case of the Navier Stokes equation written in Eulerian coordinates. As we see in Section \ref{sus:hyp_NS}, the analysis developed in the present paper applies to the study of the viscous profiles of the Navier Stokes.  

Viscous profiles provide useful information when studying the parabolic approximation of an hyperbolic conservation law. If one restricts to systems with small enough total variation, it is often meaningful to focus on viscous profiles lying  on a center, a center stable or on the uniformly stable manifold. The literature concerning these topics is extremely wide. Here, we just refer to the books by Dafermos \cite{Daf:book} and by Serre \cite{Serre:book} and to the rich bibliography contained therein for a discussion about the parabolic approximation of conservation laws. For the applications of the viscous profiles to the study of this approximation, see for example Bianchini and Bressan \cite{BiaBrevv} and Ancona and Bianchini \cite{AnBia}. Concerning the analysis of viscous profiles, we only refer to Benzoni-Gavage, Rousset, Serre and Zumbrun \cite{BenRouSerreZum}, to Zumbrun \cite{Zum:review}, and to the references therein. 
For an alternative approach to the analysis of the viscous profiles of the compressible Navier Stokes equation,  see Wagner \cite{Wagner} and the references therein. 
\vspace{1cm}

The exposition is organized as follows. In Section \ref{sus:intro_toy} we discuss a toy model and we introduce the results that are extended in the following sections to the general case. 

In Section \ref{s:hyp} we define our hypotheses and in Section \ref{sus:hyp_NS} we show that they are satisfied by the viscous profiles of the compressible Navier Stokes equation. Also, in Section \ref{sus:hyp_count} we discuss two examples: each of them show that, if one different hypothesis is not satisfied, then the first derivative $d U/ dt$ of a solution $U$ of \eqref{e:intro:theq} may blow up in finite time.

In Section \ref{s:us} we define a class of invariant manifolds for an equation with no singularity in it (Theorem \ref{p:unif} and Proposition \ref{p:exp:d}). This analysis is applied in Section \ref{s:invman} to study the singular ODE \eqref{e:intro:theq}. In particular, in Section \ref{sus:slowfast} we define the notions of  slow and fast dynamic  and we extend the definition of center manifold to the case of the singular ODE \eqref{e:intro:theq}. In Section \ref{s:appl} we discuss how to extend the notions of uniformly stable and center stable manifold: the main result here is Theorem \ref{t:app:man}.  Finally, Section \ref{sus:m:change} is devoted to the proof of Proposition \ref{p:change}, a technical result which reduces our system to a more convenient form.

\subsection{A toy model}
\label{sus:intro_toy}
In this section we discuss a toy model for system \eqref{e:intro:theq}. The goal is introducing in a simplified context the analysis that is extended in Section \ref{s:invman} to the general case. All the conditions introduced in Section \ref{s:hyp} are satisfied by the toy model discussed here, except for Hypothesis \ref{h:cutoff}. Hypothesis \ref{h:cutoff} is a technical condition and prescribes that the functions $\phi_s$ and $\phi_{ns}$ in \eqref{e:intro:theq} satisfy the following: $\phi_s (U) + \zz (U) \phi_{ns}(U)$ is identically $\vec 0$ when $U$ is out of a small enough neighbourhood of the origin. 

Actually, our toy model can be handled with known geometric singular perturbation theory techniques. Indeed, in the present section we assume that the function $\zeta$ in \eqref{e:intro:theq}  is just a parameter, namely
\begin{equation}
\label{e:i:par}
         \frac{d \zeta}{ d t } = 0.
\end{equation}
Also, we focus on the case of  a linear system:
\begin{equation}
\label{e:i:lin}
          \frac{d V}{ d t } = \frac{1}{\zeta } A_s V + A_{ns } V \qquad V \in \mathbb R^d.
\end{equation} 
In the following we consider only non negative values of $t$ and we focus on the limit $\zeta \to 0^+$. The study of the limit $\zeta \to 0^-$ does not involve additional difficulties. 

Consider the system 
\begin{equation}
\label{e:i:tau}
          \frac{d V}{ d \tau } =   A_s V +  \zz A_{ns } V \qquad V \in \mathbb R^d,
\end{equation}
which is obtained from \eqref{e:i:lin} \emph{via} the change of variable $\tau = t \ \zz$. In the following, we denote by $n_-$ the number of eigenvalues of $A_s$ having strictly negative real part (each of them counted according to its multiplicity) and by $n_+$ the number of eigenvalues with strictly positive real part. We denote by $n_0$ the multiplicity of the eigenvalue $0$ and, relying on Assumption 1 in the introduction, we assume that there are no purely imaginary eigenvalues. Also, if we write the Jordan form of $A_s$, then in the block corresponding to the eigenvalue $0$ all the entries are $0$.

Let  $\zz \to 0^+$: we are concerned with the behavior of the eigenvalues of the matrix  $A_s + \zz A_{ns}$. Thanks to results concerning the perturbation of finite-dimensional linear operators (see for example the book by Kato \cite{Kato}, page 64 and followings), these eigenvalues can be classified as follows: 
\begin{enumerate}
\item  $n_-$ eigenvalues converge to the eigenvalues of $A_s$ with strictly negative real part. We denote by $M^-(\zz)$ the eigenspace of $A_s + \zz A_{ns}$ associated to these eigenvalues. \\
\item  $n_+$ eigenvalues converge to the eigenvalues of $A_s$ with strictly positive real part. We denote by $M^+(\zz)$ the eigenspace of $A_s + \zz A_{ns}$ associated to these eigenvalues.\\
\item the remaining $n_0$ eigenvalues converge to $0$ as $\zz \to 0^+$. We denote by $M^0(\zz)$ the eigenspace of $A_s + \zz A_{ns}$ associated to these eigenvalues.
\end{enumerate} 
When $\zz \to 0^+$, the subspace $M^- (\zz)$ converges to $M^-(0)$, which is the eigenspace of $A_s$ associated to eigenvalues with strictly negative real part. The convergence occurs in the following sense: $M^-(\zz)$ is the range of a linear application $P^-(\zz) \in \mathcal L (\mathbb R^d, \; \mathbb R^d)$. As $\zz \to 0^+$, $P^-(\zz)$ converges to $P^-(0)$ and the range of $P^-(0)$ is exactly $M^-(0)$. Similarly,  when $\zz \to 0^+$, the subspaces $M^+ (\zz)$ and $M^0(\zz)$ converge respectively to $M^+(0)$ and $M^0(0)$, the eigenspaces of $A_s$ associated to the eigenvalues 
with strictly positive and zero real part. We refer again to Kato \cite{Kato} for a complete discussion. 

If $V$ belongs to $M^- (\zz)$, then 
\eqref{e:i:tau} is equivalent to 
$$
      \frac{d V^-}{ d \tau \, } =   \Big[ A^-_s  +  \zz \, \unpo \Big] V^-,
$$
where $V^- \in \mathbb R^{n_-}$ and $A^-_s$ is a $n_- \times n_-$-dimensional 
matrix which does not depend on $\zz$ and whose eigenvalues have all strictly negative real part. In the previous equation, the entries of the vector $V^-$ are the coordinates of $u$ with respect to a basis of 
$M^-(\zz)$ and $\unpo$ denotes a $n_- \times n_-$-dimensional 
matrix which possibly depends on $\zz$ but remains bounded as $\zz \to 0^+$. Its exact expression is not important here.  

If $\zz$ is sufficiently small, then all the eigenvalues of the matrix $\big[ A_s + \zz \unpo \big]$ have strictly negative real part and hence the solution $V^-(\tau)$ converges exponentially fast to $\vec 0$. More precisely, one has 
$$
   |V^- (\tau) | \leq e^{- c \tau/2} |V^-(0)|,
$$ 
where $c >0$ satisfies 
$
    - c > \lambda
$
for every $\lambda$ eigenvalue of $A_s$.  Coming back to the original variable $t$, $V^-$ satisfies 
$$
   |V^- (t) | \leq e^{- c t /2 \zz} |V^-(0)| .
$$ 
and hence the speed of exponential decay gets faster and faster as $\zz \to 0^+$. In this sense, we can regard $V^-$ as a \emph{fast dynamic}.

Conversely, assume that $V \in M^0 (\zz)$, then \eqref{e:i:tau} is equivalent to
$$
    \frac{d V^0}{ d \tau} = \zz  \Big[ L_0 A_{ns} R_0  +  \zz \unpo \Big] V^0, 
$$
where $R_0$ and $L_0$ are two matrices that do not depend on $\zz$. The matrix $R_0$ has dimension $N \times n_0$ and its columns constitute a basis of $M^0(0)$. The matrix $L_0$ is $n_0 \times N$-dimensional and satisfies
$L_0 R_0 = I_{n_0}$. Also, $V^0 = L^0 V$ and $\unpo$ denotes an $n_0 \times n_0$-dimensional matrix, which possibly depends on $\zz$ but remains bounded as $\zz \to 0^+$. Its exact expression is not relevant here. Coming back to the original variable $t$, one gets 
\begin{equation}
\label{e:i:nonos}
    \frac{d V^0}{ d t} =  \Big[ L_0 A_{ns} R_0  +  \zz \unpo \Big] V^0,
\end{equation}
and hence $V^0$ can be regarded as a \emph{slow dynamic}, because it satisfies the non singular ODE \eqref{e:i:nonos}. 

Applying the same techniques mentioned before, one gets that the eigenvalues of $ L_0 A_{ns} R_0  +  \zz \unpo$ can be divided into 3 groups:
\begin{enumerate}
\item eigenvalues that converge to the eigenvalues of  $L_0 A_{ns} R_0$ with strictly negative real part. We denote by $M^{0-}(\zz)$ the corresponding eigenspace.
\item  eigenvalues that converge to the eigenvalues of  $L_0 A_{ns} R_0$ with strictly positive real part. We denote by $M^{0+}(\zz)$ the corresponding eigenspace. 
\item  eigenvalues that converge to the eigenvalues of  $L_0 A_{ns} R_0$ with zero part. We denote by $M^{00}(\zz)$ the corresponding eigenspace. 
\end{enumerate}
If $V (t)  \in M^{0-} (\zz) $, then $V(t)$ converges exponentially fast to the equilibrium $\vec 0$ when $ t \to + \infty$, but the speed of exponential decay does not blow up as $\zz \to 0^+$. 

The space 
\begin{equation}
\label{e:i:sum}
         M^s (\zz) = M^{-} (\zz)  \oplus M^{0-} (\zz)
\end{equation}      
can be regarded as \emph{uniformly stable space} for \eqref{e:i:lin}, because every orbit entirely contained in this space decays exponentially fast to $\vec 0$. Also, the speed of exponential decay is uniformly bounded from below by a constant which does not depend on $\zz$. Another way of interpreting this observation is the following: combine equations \eqref{e:i:par} and \eqref{e:i:lin} and consider in $\mathbb R^{d +1}$ the system 
\begin{equation}
\label{e:i:dplus}
\left\{
\begin{array}{ll}
            d \zz / dt = 0 \\
            d V / dt = A_s V / \zz + A_{ns} V.  \\
\end{array}
\right.
\end{equation}
Every point in the subspace $\{ (\zz, \, \vec 0) \}$ is then an equilibrium for \eqref{e:i:dplus}. Also, for any $\bar \zz$, every orbit $V(t)$ lying on $M(\bar \zz)$ converges to the equilibrium $(\bar \zz, \, \vec 0)$, and the speed of exponential decay is bounded below by a constant independent of $\bar \zz$.

Conversely, the space $M^{00} (\zz)$ can be regarded as a \emph{center manifold} of the original equation \eqref{e:i:lin}.

In Section \ref{s:invman} we will extend the previous considerations to the case of the non linear equation 
$$
     \frac{d U}{d t} = \frac{1}{\zeta (U)} \phi_s (U) + \phi_{ns} (U),
$$
where $\zz (U)$ is in general a non constant function. 
\section{Hypotheses}
\label{s:hyp}
In this section we discuss the hypotheses assumed in the work. 

More precisely, in Section \ref{sus:hyp_hyp} we state the conditions imposed on the singular ODE
\begin{equation}
\label{e:hyp:sin}
          \frac{d U}{d t} = \frac{1}{\zeta (U)} \phi_s (U) + \phi_{ns} (U).
\end{equation}
These conditions can be divided into two groups: Hypotheses \ref{h:pos}, \ref{h:cutoff}, \ref{h:neg}, allow to avoid some technical complications, but could be  actually omitted at the price of much heavier notations. On the other side, Hypotheses \ref{h:sur}, \ref{h:center}, \ref{h:tras}, \ref{h:fast}, \ref{h:slow} are much more important and they will be deeply exploited  in Section \ref{s:invman}. Note, however, that in Section \ref{s:us} we are not directly concerned with the singular ODE \eqref{e:hyp:sin} and that we do not exploit Hypotheses  \ref{h:sur}, \ref{h:center}, \ref{h:tras}, \ref{h:fast}, \ref{h:slow}.

Moreover, in Section \ref{sus:hyp_count} we discuss three counterexamples. They show that, if either Hypothesis \ref{h:fast} or Hypothesis \ref{h:slow} is violated, then the results discussed in the following sections do not hold. In particular, there might be solutions of \eqref{e:hyp:sin} that are not continuously differentiable.

Finally, in Section \ref{sus:hyp_NS} we verify that the conditions introduced in Section \ref{sus:hyp_hyp} are satisfied by the viscous profiles of the compressible  Navier Stokes equation written in Eulerian coordinates. 

\subsection{Hypotheses satisfied by the singular O.D.E}
\label{sus:hyp_hyp}
Define 
\begin{equation}
\label{e:hyp:F}
           F (U) =  \phi_s (U) + \zeta (U) \phi_{ns} (U) 
\end{equation}
and consider the non singular ordinary differential equation 
\begin{equation}
\label{e:hyp:theq}
           \frac{d U}{d \tau} = F(U)  \qquad U \in \mathbb R^N.
\end{equation}
Formally, \eqref{e:hyp:theq} is obtained from \eqref{e:hyp:sin} {\it via} the change of variables $\tau = \tau (t)$ defined as the solution of the Cauchy problem
\begin{equation}
\label{e:hyp:cauchy}
\left\{
\begin{array}{lll}
          \displaystyle{   \frac{d \tau}{dt} = \frac{1}{\zeta [U(t) ]} } \\
         \\
          \tau ( 0) =0.\\
\end{array}
\right.
\end{equation}
However, the function $\tau (t)$ is well defined only if $\zeta [U(t)] \neq 0 $ for every $t$. To overcome this difficulty we then proceed as follows: we state all the hypotheses referring to the formulation \eqref{e:hyp:theq}. Relying on these hypotheses, in Section \ref{s:invman} we prove the existence of various locally invariant manifolds for \eqref{e:hyp:theq} satisfying the following property. If $U$ is an orbit lying on one of these manifolds and $\zeta [U(0)] \neq 0$, then $\zeta [U(t)] \neq 0$ for every $t$. If we restrict to the orbits lying on these manifolds, \eqref{e:hyp:theq} is equivalent to \eqref{e:hyp:sin}.  

To simplify the exposition, we assume the following:
\begin{hyp}
\label{h:pos}
         The initial datum $U(0)$ of \eqref{e:hyp:theq} satisfies $\zz \big( U (0) \big) > 0$.  
\end{hyp}
The case $\zz \big( U(0) \big) < 0$ does not involve additional difficulties. The main difference is that, if  $\zz \big( U(0) \big) < 0$, then the change of variable defined by \eqref{e:hyp:cauchy} has negative derivative. As a consequence, when $t \to + \infty$ the function $\tau (t) \to - \infty$. Loosely speaking, the statements given in the present paper can be extended to the case $\zz \big( U(0) \big) < 0$ in the following way. All the statements concerning the fast dynamics and referring to the \emph{stable} space or to \emph{stable}-like manifolds have to be replaced by analogous statements concerning the \emph{unstable} space or \emph{unstable}-like manifolds. However, we will not consider the case $\zz \big( U(0) \big) < 0$ explicitly. 

Before stating the other hypotheses, we recall that we want to study \eqref{e:hyp:sin} and \eqref{e:hyp:theq} in the neighbourhood of an equilibrium point $\bar U$ such that 
$F(\bar U) = \vec 0$ and $\zeta (\bar U)= 0$. It is not restrictive to take $\bar U = \vec 0$. Namely, in the following we assume
\begin{equation}
\label{e:hyp:0}
         F (\vec 0) = \vec 0 \qquad \zeta (\vec 0) = 0 . 
\end{equation}
Also, we can assume the following. Fix a positive constant $\delta>0$ and consider a smooth cut-off function $\rho (U)$ satisfying 
$$
    \rho (U) =
    \left\{
    \begin{array}{ll}
               U & |U| \leq \delta \\
               0 &  |U| \ge 2 \delta.
    \end{array}
    \right.
$$
In the following, instead of studying system \eqref{e:hyp:theq} we focus on 
$$
     \frac{d U}{ d \tau} = \rho (U) F(U).
$$
However, to simplify the notations instead of writing each time $ \rho (U) F(U)$ we assume that Hypothesis 
\ref{h:cutoff} holds. 
\begin{hyp}
\label{h:cutoff}
             The function $F$ satisfies the following condition: if $|U| \ge 2 \delta$ then $F(U) = \vec 0$.
\end{hyp}
The exact size of the constant $\delta$ will be discussed in the following.

Note that Hypothesis \ref{h:cutoff} is not restrictive if the goal is to study the solutions of \eqref{e:hyp:theq} that remain confined in a neighbourhood of the origin of
 size $\delta$. Loosely speaking, the analysis developed in Sections \ref{s:us} and \ref{s:invman} can be extended to the orbits of systems that violate Hypothesis \ref{h:cutoff} as far as these orbits remain in a neighbourhood of the origin with size $\delta$. In particular, the manifold described in Sections \ref{s:us} and \ref{s:invman} are no more \emph{invariant} if Hypothesis \ref{h:cutoff} is violated: they are just \emph{locally invariant}. 

We also introduce the following simplification. It is not restrictive to assume that all the eigenvalues of $D F (\vec 0 )$ have non positive real part. Indeed, this condition is satisfied if we restrict to a \emph{center-stable} manifold for \eqref{e:hyp:theq}.  As mentioned in the introduction, the existence of a center stable manifold can be obtained as a consequence of the Hadamard Perron theorem, which is discussed for example in the book by Katok and Hasselblatt \cite{KHass} (Chapter 6, page 242).  
Also, note that if $\zz \big( U (0) \big) <0$ then it is not restrictive to assume that all the eigenvalues of $D F (\vec 0 )$ have \emph{non negative} real part: this can be obtained considering the solutions that lie on a \emph{center unstable} manifold.

\begin{hyp}
\label{h:neg}
         The Jacobian $DF (\vec 0)$ admits only eigenvalues with non positive real part. 
\end{hyp}

Also, we assume the following non degeneracy condition:
\begin{hyp}
\label{h:sur}
         The gradient $\nabla \zeta (\vec 0)  \neq \vec 0$. 
\end{hyp}
Let $\mathcal S$ be the singular set
\begin{equation}
\label{e:hyp:S}
         \mathcal S : = \big\{ U: \; \zeta (U) =0   \big\}.
\end{equation}
Thanks to the implicit function theorem, Hypothesis \ref{h:sur} ensures that in a small enough neighbourhood of $\vec 0$ the set $\mathcal S$ is actually a manifold of dimension $N-1$, where $N$ is the dimension of $U$.  
\begin{hyp}
\label{h:center}
          Let $\mathcal M^c$ be any center manifold for \eqref{e:hyp:theq} around the equilibrium point $\vec 0$. If $|U| \leq \delta $ and $U$ belongs to the intersection $\mathcal M^c \cap \mathcal S$ , then $U$ is an equilibrium for  \eqref{e:hyp:theq}, namely $F(U)= \vec 0$ .   
\end{hyp}
Concerning equilibria, we also assume the following
\begin{hyp}
\label{h:tras}
         There exists  a manifold of equilibria  $\mathcal M^{eq}$ for \eqref{e:hyp:theq} which contains $\vec 0$ and which is transversal to $\mathcal S$.
\end{hyp}
Let $n_{eq}$ be the dimension of $\mathcal M^{eq}$. We recall that the manifolds $\mathcal S$ and $\mathcal M^{eq}$ are transversal if the intersection $\mathcal S \cap \mathcal M^{eq}$ is a manifold with dimension $n_{eq} - 1$ (as pointed out before, the dimension of $\mathcal S$ is $N-1$).   
\begin{hyp}
\label{h:fast}
          For every $U \in \mathcal S$, 
          \begin{equation}
          \label{e:hyp:fast}
                    \nabla \zeta (U) \cdot F(U) =0  .   
          \end{equation}
\end{hyp}
Thanks to Hypothesis \ref{h:fast} and to the regularity of the functions $\zeta$ and $F$, the function
$$
    G(U) = \frac{   \nabla \zeta (U) \cdot F(U)}{ \zeta (U)}
$$
can be extended and defined by continuity on the surface $\mathcal S$. 
\begin{hyp}
\label{h:slow}
          Let $U \in \mathcal S$ be an equilibrium for \eqref{e:hyp:theq}, namely $\zeta (U)= 0$ and $F(U) = \vec 0$. Then 
          \begin{equation}
          \label{e:hyp:slow}
                   G(U) = 0. 
          \end{equation}
\end{hyp}
In Section \ref{sus:intro_toy} we introduced, in the case of a toy model, the notion of \emph{slow} and \emph{fast} dynamics. These notions will be extended in Section \ref{sus:slowfast} to the general non linear case. Hypotheses \ref{h:slow} and \ref{h:fast} can be then reformulated saying that the set $\mathcal S$ is invariant for the manifold of the slow and of  the fast dynamics respectively.  
\begin{rem}
         Consider system
         \begin{equation}
         \label{e:hyp:r:old}
                  \frac{d U}{d t } = \frac{1}{ \zeta (U)} \phi_s (U) + \phi_{ns} (U). 
         \end{equation}
         Also, let $f(U)$ be a regular, real valued function such that $f(\vec 0) > 0$. Clearly, \eqref{e:hyp:r:old} is equivalent to 
           \begin{equation}
         \label{e:hyp:r:new}
                  \frac{d U}{d t } = \frac{1}{ \zeta (U) f(U) } \phi_s (U) f (U) + \phi_{ns} (U)
         \end{equation}
         and $\zeta (U) f (U) \to 0^+$ if and only if $\zeta (U) \to 0^+$, at least in a sufficiently small neighbourhood of $U = \vec 0$.  By direct check, one can verify that Hypotheses \ref{h:pos} $\dots$ \ref{h:slow} are verified by the couple $( \zeta, \, F )$ if and only if they are verified by the couple 
         $( \zeta f , \,  F f )$. 
\end{rem}

\begin{rem} 
	As we will see in Section \ref{s:invman}, Hypothesis \ref{h:center} can be reformulated saying that the slow dynamics intersecting the singular manifold $\{ U: \; \zz (U) = 0 \}$ are equilibria for system \eqref{e:hyp:theq}.  Heuristically, this means that we require that the limit as $\zz \big( U(0) \big) \to 0^+$ of a solution of \eqref{e:hyp:sin} is a solution of the limit system. In other words, we want to rule out the possibility of a relaxation effect. 

	As it shown by the examples discussed in next section, the assumptions on the invariance of the manifold  $\{ U: \; \zz (U) = 0 \}$ with respect to the slow and the fast dynamics (Hypothesis \ref{h:fast} and \ref{h:slow} respectively) are due to the fact that we want to have a smooth invertible time rescaling $t= t(\tau)$ defined by \eqref{e:hyp:cauchy}.
\end{rem}

\subsubsection{The case of the compressible Navier Stokes in Eulerian coordinates}
\label{sus:hyp_NS}
In this section we show that Hypotheses \ref{h:pos}, \ref{h:neg}, $\dots$ \ref{h:slow} are satisfied by the ODE for the viscous profiles of the compressible 
Navier Stokes equation written in Eulerian coordinates. Also, Hypothesis \ref{h:cutoff} is not restrictive if the goal is to study the viscous profiles entirely contained in a small neighbourhood of an equilibrium point. 

The case of the Navier Stokes written in Lagrangian coordinates  was already discussed 
for example in Rousset \cite{Rousset:char}. When the equation is formulated using Lagrangian coordinates, the ODE satisfied by the viscous profiles is not singular. 

The compressible Navier Stokes written in Eulerian coordinates is
\begin{equation}
\label{e:ns:eul}
       \left\{
       \begin{array}{lll}
              \rho_t + ( \rho v )_x =0 \\
	      (\rho v)_t + \Big( \rho v^2 + p \Big)_x = \displaystyle{ \Big( \nu v_x  \Big)_x } \\
	      \displaystyle{ \Big( \rho e + \rho \frac{v^2}{2}\Big)_t + \Big(v \Big[ \frac{1}{2} \rho v^2 
	      + \rho e + p \Big] \Big)_x = \Big( k \theta_x + 
	      \nu v v_x \Big)_x}. \\
       \end{array}
       \right.
\end{equation}
Here, the unknowns are $\rho(t, \, x), \, v(t, \, x)$ and $\theta (t, \, x)$. The function $\rho$ represents the density of the fluid, $v$ is the velocity of the particles in the fluid and $\theta$ is the absolute temperature. The function $p= p(\rho, \, \theta) >0$ is the pressure and satisfies $p_{\rho} >0$, while $e$ represent the internal energy. In the case of a polytropic gas, the following relation holds:
$
   e =\theta R / (\gamma -1 ),  
$    $R$ being the universal gas constant and $\gamma >1$ a constant specific of the gas. Finally, $\nu(\rho)>0$ and $k(\rho)>0$ represent the viscosity and the heat conduction coefficients respectively.

After some manipulations (see \cite{BiaSpi:pro} for details), one gets that the equation satisfied by the steady solutions of the compressible Navier Stokes can be written in the form 
$$
    \frac{d U}{ d x } = \frac{1}{ \zz (U)} F(U)
$$
provided that $U = \Big( \rho, \, v, \, \theta, \, \vec z \Big)^T$, $\zz (U) = v$ and 
\begin{equation}
\label{e:ns:f}
    F(U) =
    \left(
    \begin{array}{ccc}
                A_{21}^t \vec{z} / a_{11} \\
                v \, \vec z   \\
                 b^{-1} \displaystyle{\Big[     A_{22} v   -   A_{21} A_{21}^T  /a_{11}   \Big] \vec z } \\
    \end{array}
    \right)
\end{equation}
The equation satisfied by the travelling waves of the compressible Navier Stokes equation in one space variable is similar, the only difference being that the singular value is $v= \sigma$, where $\sigma$ is the speed of the travelling wave. 

In \eqref{e:ns:f}, $A_{21}$ is a vector in $\mathbb R^2$ and $A_{21}^t$ denotes its transpose. Also, $ w = \rho _x$  ad $\vec z = \big( v_x, \, \theta_x \big)^t$. The function $a_{11}$ is real valued and strictly positive if $\rho$ is bounded away from $0$.  The matrix $b$ has dimension $2 \times 2$ and all its eigenvalues have strictly positive real part. The exact expression of these terms is not important here. Finally,  
\begin{equation}
\label{e:hyp:a22}
    A_{22} = \frac{1}{ \theta }
     \left(
    \begin{array}{ccc}
              \rho v -  \nu' \rho_x &  p_{\theta} \\
                      p_{\theta}  - 
                       \nu v_x /  \theta  
                      &    \rho v e_{\theta} / \theta  
                       -  k' \rho_x  /  \theta  \\              \end{array}
    \right)  
\end{equation}
Note that $A_{22}$ depends on $\rho_x$ but, plugging $  w =  - A_{21}^T \vec{z} / ( a_{11}   v )$ into \eqref{e:hyp:a22} one gets that $A_{22}v$ evaluated at a point $(\rho, \, v=0, \, \theta, \, \vec z = \vec 0)$
is the null matrix. 

Thus, the Jacobian $DF$ satisfies 
$$
     DF (\rho, \, v=0, \, \theta, \, \vec z = \vec 0) =
     \left(
     \begin{array}{ccccc}
                 0 & -  A_{21}^T   / a_{11}   \\
                 \\
                 \vec 0 &  \displaystyle{ \mathbf{0}}_{2} \\
                 \\
                 \vec 0 &    - b^{-1}    A_{21} A_{21}^T  /a_{11}    \\
     \end{array}
     \right),
     $$
where $ \displaystyle{ \mathbf{0}}_{2}$ denotes the $2 \times 2$ null matrix. Since 
$  A_{21} A_{21}^T  /a_{11} $ admits only eigenvalues with strictly positive real part, then $DF$ admits only eigenvalues with non positive real part and hence Hypothesis \ref{h:neg}
is satisfied. Also, the dimension of every center manifold of 
$$
      \frac{d U}{ d \tau } = F(U)
$$
is $3$. Since the subspace $\{ \vec z = \vec 0 \}$ is entirely constituted by equilibria of the equation, it coincides with the center manifold. Thus, Hypothesis \ref{h:center} is satisfied. 
Since $\zz (U) = v$, then Hypothesis \ref{h:sur} is also verified. Concerning Hypothesis \ref{h:tras}, this is satisfied because $\{ \vec z = \vec 0 \}$ is transversal to the singular surface 
$\{ v =0 \}$. Finally, by direct check one can show that also Hypotheses \ref{h:fast} and \ref{h:slow} are verified. Thus the machinery developed in this work applies to the study of viscous profiles with small total variation of the compressible Navier Stokes equation written in Eulerian coordinates.

\subsection{Examples}
\label{sus:hyp_count}

\subsubsection{Example \eqref{e:ex:fast}}
\label{e:fast}
          Example \eqref{e:ex:fast} deals with a system which satisfies Hypotheses \ref{h:pos}, \ref{h:neg} $\dots$ \ref{h:tras} and Hypothesis \ref{h:slow}, but does not satisfy 
Hypothesis \ref{h:fast}. We exhibit a solution of this system which has a blow up in the first derivative and hence it is not continuously differentiable. The loss of regularity experienced in Example 
\eqref{e:ex:fast} regards a solution $U$ such that $\zeta [U (0)]\neq 0$, but $\zeta (U)$ reaches the value $0$ for a finite value of $t$.

           Consider the following system:
           \begin{equation}
           \label{e:ex:fast}
           \left\{
           \begin{array}{ll}
                      d u_1 / dt = - u_2 / u_1 \\
                      d u_2 / dt = - u_2 \\
           \end{array}
           \right.           
           \end{equation}         
Let $U = \Big( u_1, \, u_2 \Big)^T$, $\zeta (U) = u_1$ and 
$$
    \phi_s (U) =
    \left(
    \begin{array}{cc}
                - u_2 \\
                0 \\
     \end{array}
     \right)
     \qquad 
      \phi_{ns} (U) =
    \left(
    \begin{array}{cc}
                0  \\
                - u_2  \\
     \end{array}
     \right).               
$$
System \eqref{e:ex:fast} can then be written in the form 
$$
     \frac{d U}{d t} = \frac{1}{\zeta (U)} \phi_s (U) + \phi_{ns} (U).
$$     
In this case, the function $F(U)$ defined by  \eqref{e:hyp:F}
is 
$$
    F(U) = 
    \left(
    \begin{array}{cc}
                - u_2   \\
                - u_2 u_1 \\
     \end{array}
     \right). 
$$
By direct check, one can verify that Hypotheses  \ref{h:pos} $\dots$ \ref{h:tras} and Hypothesis \ref{h:slow} are satisfied by \eqref{e:ex:fast}. On the other side, Hypothesis \ref{h:fast}
is not verified in this case. Indeed, the singular surface $\mathcal S$ defined by \eqref{e:hyp:S} is in this case the line $\{ u_1 = 0 \}$ and 
$$
    \nabla \zeta \cdot F = - u_2 
$$
is in general different from $0$ on $\mathcal S$. 

The solution of \eqref{e:ex:fast} can be explicitly computed and it is given by 
\begin{equation}
\label{e:ex:sol}
     \left\{
    \begin{array}{lll}
                \displaystyle{ u_1 (t) = \sqrt{  u_1 (0) +  u_2 (0) \big(  e^{- t } - 1 \big)  }} \\
                \\
                   \displaystyle{  u_2 (t) = u_2 (0) e^{- t } } \\
    \end{array}
    \right.
\end{equation}
Choosing $u_2 (0) > u_1 (0) >0$, one has that the solution $u_1 (t)$ can reach the value $0$ for a finite $t$. Note that at that point $t$ the first derivative $d u_1 / dt$
blows up: thus, the solution \eqref{e:ex:sol} of \eqref{e:ex:fast} is not $\mathcal C^1$. 
\subsubsection{Example \eqref{e:ex:ok}}
\label{e:ok}
          Example \eqref{e:ok} deals with system \eqref{e:ex:ok}, which is apparently very similar to \eqref{e:ex:fast}. However, in the case of \eqref{e:ex:ok} Hypotheses \ref{h:pos}, \ref{h:neg} $\dots$ \ref{h:slow} are all verified. We show the solutions of \eqref{e:ex:ok} are regular. Also, if $\zeta [ U(0) ] \neq 0$ then $\zeta [U(t) ] \neq 0$ for all values of $t$.

          Consider system 
          \begin{equation}
          \label{e:ex:ok}
                     \left\{
                     \begin{array}{cc}
                                 d u_1 / dt = - u_2 \\
                                 d u_2  / dt = - u_2 / u_1
                     \end{array}
                     \right.
          \end{equation}
          Set $U = \Big( u_1, \, u_2 \Big)^T$, $\zeta (U) = u_1$.  System \eqref{e:ex:ok} can be written in the form 
$$
     \frac{d U}{d t} = \frac{1}{\zeta (U)} \phi_s (U) + \phi_{ns} (U)
$$   
provided that 
$$
    \phi_s (U) =
    \left(
    \begin{array}{cc}
                0 \\
                - u_2  \\
     \end{array}
     \right)
     \qquad 
      \phi_{ns} (U) =
    \left(
    \begin{array}{cc}
                - u_2  \\
                  0  \\
     \end{array}
     \right).               
$$  
Also, the function $F(U)$ defined by  \eqref{e:hyp:F}
is in this case 
$$
    F(U) = 
    \left(
    \begin{array}{cc}
                - u_2 u_1    \\
                - u_2 \\
     \end{array}
     \right). 
$$
By direct check, one can verify that Hypotheses \ref{h:pos} $\dots$ \ref{h:slow} are all verified in this case. 

To study system \eqref{e:ex:ok} we can proceed as follows. From \eqref{e:ex:ok} we have 
$$
    \frac{d u_1 / dt }{ u_1 } = - \frac{u_2 }{ u_1 } = d u_2 / dt 
$$
and hence 
$$
     \ln \Bigg[ \frac{u_1 (t)}{u_1 (0)} \Bigg] = u_2 (t) - u_2(0).
$$
Eventually, we obtain
\begin{equation}
\label{e:ex:u1}
    \displaystyle{ u_1 (t) = u_1 (0) e^{u_2 (t) - u_2 (0) } }.
\end{equation}
Choose $u_1 (0) > 0$. To prove that $u_1 (t) \neq 0$ for all $t$ it is enough to show that $u_2 (t)$ is well defined (and in particular finite) for every $t>0$. In the following we also prove 
that $u_2 (t)$ is also $\mathcal C^{\infty}$ for every $ t \ge 0$. This guarantees that no loss of regularity occurs.

Plugging \eqref{e:ex:u1} into the second line of \eqref{e:ex:ok} we get
\begin{equation}
\label{e:ex:u2}
     \displaystyle{     d u_2 / dt  = - \frac{u_2 }{u_1 (0)} e^{u_2 (0) - u_2 (t) } }. 
\end{equation}
Note that $u_2 = 0$ is an equilibrium for \eqref{e:ex:u2}. Also, if $u_2 (0) < 0$ then $d u_2 / dt \ge 0$ and hence $u_2 (0) \leq u_2 (t) < 0$ for every $t$. Conversely, if 
$u_2 ( 0) > 0$ then $d u_2 / dt \leq  0$ and hence $0 \leq u_2 (t) < u_2 (0)$ for every $t$.  In both cases, we get that $u_2 (t)$ is well defined and regular for every $t \ge 0$. 
\subsubsection{Example \eqref{e:ex:slow}}
\label{e:slow}
            With Example \eqref{e:slow} we discuss a system which satisfies 
Hypotheses \ref{h:pos}, \ref{h:neg} $\dots$ \ref{h:fast}, but does not satisfy 
Hypothesis \ref{h:slow}. As in Example \eqref{e:ex:fast}, we exhibit a solution of this system which is not continuously differentiable and the loss of regularity regards a solution $U$ such that $\zeta [U (0)]\neq 0$, but $\zeta (U)$ reaches the value $0$ for a finite value of $t$. 
 
           Consider the following system:
           \begin{equation}
           \label{e:ex:slow}
           \left\{
           \begin{array}{lll}
                      d u_1 / dt = - u_3 \\
                      d u_2 / dt = - u_2 / u_1  \\
                      d u_3 / dt = - u_3 \\
           \end{array}
           \right.           
           \end{equation}         
Let $U = \Big( u_1, \, u_2, \, u_3  \Big)^T$, $\zeta (U) = u_1$ and 
$$
    \phi_s (U) =
    \left(
    \begin{array}{ccc}
                  0  \\
                - u_2  \\
                 0 \\
     \end{array}
     \right)
     \qquad 
      \phi_{ns} (U) =
    \left(
    \begin{array}{ccc}
                  - u_3  \\
                0   \\
                - u_3   \\  
      \end{array}
     \right).               
$$
System \eqref{e:ex:fast} can then be written in the form 
$$
     \frac{d U}{d t} = \frac{1}{\zeta (U)} \phi_s (U) + \phi_{ns} (U)
$$     
and the function $F(U)$ defined by  \eqref{e:hyp:F}
is 
$$
    F(U) = 
    \left(
    \begin{array}{ccc}
                - u_3 u_1   \\
                - u_2  \\
                - u_3 u_1
     \end{array}
     \right). 
$$
By direct check, one can verify that Hypotheses  \ref{h:pos} $\dots$ \ref{h:fast} are verified by \eqref{e:ex:slow}. On the other side, Hypothesis \ref{h:slow}
is not satisfied in this case. Indeed, the surface $\mathcal S = \{ U: \; \zeta (U) =0 \}$ is the plane $\{ u_1 = 0 \}$. Thus, the set of points such that $\zeta (U) =0$ and 
$F (U) = \vec 0$ is $\{ u_1 = u_2 = 0 \}$ and 
$$
    \nabla \zeta \cdot DF \cdot  \Big( \nabla \zeta \Big)^T = -  u_3 
$$  
is in general different from zero on this line. 

An explicit solution of \eqref{e:ex:slow} can be obtained as follows. From the third and the first equation we get respectively
\begin{eqnarray*}
         &  \displaystyle{ u_3 (t) = u_3 (0) e^{- t }} \\
         &  \displaystyle{ u_1 (t) = u_1 (0) - u_3 (0) + u_3 (0) e^{-t}  }.  \\
\end{eqnarray*}       
Assume that $u_3 (0) = A  u_1(0)$ for some constant $A$ whose exact value is determined in the following. The equation satisfied by $u_2$ becomes 
$$
    \frac{d u_2 }{d t }= - \frac{u_2 }{A  u_1 (0)  e^{- t } + u_1(0) (1- A) }. 
$$
Thus, we obtain
$$
    \frac{ d}{dt}  \Big[ \ln \Big( u_2 (t) \Big) \Big] = \frac{1}{u_1 (0) (A -1)}  \frac{ d}{ d t } \Big[    \ln \Big(  u_1 (0)(1-A)  e^{t} +A  u_1 (0) \Big)  \Big]  
$$
and hence 
$$
    u_2 (t) = B \Big[(1-A) e^t + A \Big]^{ \displaystyle{  \frac{1 }{  (A-1) u_1 (0) }}}
$$
for a suitable constant $B$. 
If $(A-1)u_1 (0) >1$, then the first derivative $d u_2 / dt$ blows up at $t = \ln (A / A-1)$. Note that this is exactly the value of $t$ at which $u_1 (t)$ attains $0$.  

In general, for every $u_1 (0) > 0$ if $1/ (A-1)u_1 (0) $ is not a natural number, then the solution is not in $\mathcal{C}^m $ for $m= [1 / (A-1)u_1 (0)  ] +1 $. Here $ [1 / (A-1)u_1 (0) ] $ denotes the entire part.  Thus,  we have a loss of regularity in higher derivatives. 
\section{Uniformly stable manifolds}
\label{s:us}
In this section we consider the system 
\begin{equation}
\label{e:us:theq}
           \frac{d U}{d \tau} = F(U),
\end{equation}
where $U \in \mathbb R^N$ and the $F: \mathbb R^N \to \mathbb R^N$ is a regular function. We are interested in the behavior of the solutions in a small enough neighbourhood of an equilibrium point.  We can then assume that such an equilibrium point is $\vec 0$, namely $F (\vec 0) = \vec 0$. Also, we assume 
\begin{equation}
\label{e:us:cutoff}
     |U| \ge 2 \delta  \implies F(U) = \vec 0. 
\end{equation}
Because of Hypothesis \ref{h:cutoff}, \eqref{e:us:cutoff} is not restrictive in view of the applications discussed in Section \ref{s:invman}. Also, because of Hypothesis \ref{h:neg} we assume that all the eigenvalues of the Jacobian $DF (\vec 0)$ have non positive real part. Note, however, that in this Section we exploit none among Hypotheses \ref{h:pos}, \ref{h:sur}, \ref{h:center}, \ref{h:tras}, \ref{h:fast} and \ref{h:slow}. 
\subsection{Notations and preliminary results}
\label{s:us:prel}
\subsubsection{Fr\'{e}chet differentiability of the fixed point of a family of maps}
\label{s:us:frechet}
In the following we have to study the regularity of the fixed points of a family of maps depending on a parameter. To do this, we exploit Lemma \ref{l:fre}. Note that the hypotheses there are not sharp and the result could be improved. However, to avoid technical complications we restrict to those hypotheses since they are satisfied in the cases we discuss in the following.  

Let $X$ be a closed subset with non empty interior in a Banach space $\tilde X$ and let $Y$ be an open subset of another Banach space $\tilde Y$. Also, let 
$$
    T: X \times Y \to \tilde X
$$
be a map such that, for every $y \in Y$, $T(\cdot, \, y)$ takes values in $X$ and is a strict contraction, namely there exists some constant $k < 1$ such that
$$
    \| T (x_1, \, y) - T(x_2, \, y)  \|_{\tilde X} \leq k \| x_1 - x_2 \|_{\tilde X} \quad \forall \, x_1, \; x_2 \in X. 
$$
Thanks to the Contraction Mapping Theorem, we can define a map 
$$
    Y \to X
$$
which associates to $y$ the fixed point of the map $T(\cdot, \, y)$. We denote this map by $x (y)$ and we are interested in its regularity. Assume that, for every $y$, $x(y)$ belongs to the inner part of $X$. Also, assume that, for every $(\bar x, \, \bar y) \in X \times Y$ such that $\bar x$ is an inner point, $T(\cdot, \, \bar y)$ is Fr\'{e}chet differentiable at $\bar x$ and denote by 
$T_x (\bar x, \, \bar y) \in \mathcal L (\tilde X, \, \tilde X)$ its differential. Also, assume that, for every  $(\bar x, \, \bar y) $ in the same conditions as before, the map $T(\bar x, \, \cdot)$ is Fr\'{e}chet differentiable at $\bar y$ and denote by $T_y (\bar x, \, \bar y) \in \mathcal L (\tilde Y, \,  \tilde X)$ its differential. 

The proof of the following result relies on standard techniques (see for example the book by Ambrosetti and Prodi \cite{AmPro}) and will be therefore omitted.  
\begin{lem}
\label{l:fre}
           Assume that the map $x (y)$ is Lipschitz continuous and fix a point $(\bar x, \, \bar y)$ such that $\bar x = x (\bar y)$.  Also, assume that $T$ is  Fr\'{e}chet differentiable with respect to the variable $y$ at the point $(\bar x, \, \bar y)$, namely 
           $T_y (\bar x, \, \bar y)$ exists. If $T_x (x, \, y)$ is defined and continuous in a neighbourhood of $(\bar x, \, \bar y)$, then $x(y)$ is 
          Fr\'{e}chet differentiable at $\bar y$ and the differential is 
          \begin{equation}
          \label{e:m:fre:fo}
              \Big[ I - T_x \big(  x (\bar y), \,  \bar y \big) \Big]^{-1} \circ T_y \big( x (\bar y), \, \bar y \big), 
          \end{equation}
          where $I$ denotes the identity. 
\end{lem}
Note that the map $ \Big[ I - T_x \big[ x (\bar y), \,  \bar y \big]  \Big]$ is invertible because $T( \cdot, \, \bar y)$ is a strict contraction on $X$. 
\begin{rem}
\label{r:lip}
          We want to give a sufficient condition to have that $x(y)$ is Lipschitz continous. Assume that there exists a constant $L$ such that, for every $y_1, \; y_2 \in Y$, 
          \begin{equation}
          \label{e:us:suff:Lip}
              \| T \big( x (y_1), \, y_1 \big)  -  T \big( x (y_1), \, y_2 \big)  \|_X \leq L \| y_1 - y_2 \|_Y. 
          \end{equation}
          Then the map $x(y)$ is Lipschitz continuous. Indeed, 
          \begin{equation*}
          \begin{split}
                       \| x(y_1) - x(y_2) \|_X &  =   \| T \big( x (y_1), \, y_1 \big)  -  T \big( x (y_2), \, y_2 \big)  \|_X  \\
           &       \leq    \| T \big( x (y_1), \, y_1 \big)  -  T \big( x (y_1), \, y_2 \big)  \|_X +
                     \| T \big( x (y_1), \, y_2  \big)  -  T \big( x (y_2), \, y_2 \big)  \|_X \\
           &           \leq   L \| y_1 - y_2 \|_Y + k   \| x(y_1) - x(y_2) \|_X \\
          \end{split}
          \end{equation*}
          Since $k <1$, we get that $x(y)$ is Lipschitz continuous.  
\end{rem}
\subsubsection{First change of variables}
\label{s:us:first:ch}
Consider system \eqref{e:us:theq}. Let $V^-$ be the eigenspace of the Jacobian $DF (\vec 0)$ associated to eigenvalues with strictly negative real part. Also, let $V^0$ be the eigenspace associated to the eigenvalues with $0$ real part. Also, fix $\mathcal V^0$, a center manifold of \eqref{e:us:theq} around the equilibrium $\vec 0$. Finally, let $\mathcal V^-$ be the stable manifold. The manifolds $\mathcal V^0$ and $\mathcal V^-$ are tangent at the origin to $V^0$ and $V^-$ respectively. Note that $\mathbb R^N = V^0 \oplus V^-$ because $DF(\vec 0)$ admits only eigenvalues with non positive real part. Thanks to the local invertibility theorem, in a sufficiently small neighbourhood of the origin we can define a local diffeomorphism $\uu$ such that the following conditions are satisfied. Let $\bar U = \uu(U)$, then $\bar U$ satisfies 
\begin{equation}
\label{e:us:theq:b}
     \frac{d \bar U}{ d \tau} = \bar f  (\bar U),
\end{equation}     
where $\bar f(\bar U) = D \uu \big( \uu^{-1} (\bar U ) \big) F \big( \uu^{-1} (\bar U) \big).$ Write $\bar U = (\bar X^-, \, \bar X^0)$, where $\bar X^0$ has the same dimension as $V^0$ and $\bar X^-$ has the same dimension as $V^-$. Then the stable manifold of \eqref{e:us:theq:b} is the subspace $\{ \bar X^0 \equiv \vec 0 \}$, while the center manifold is the subspace $\{ \bar X^- \equiv  \vec 0 \}$. In the following, we assume that the constant $\delta $ in \eqref{e:us:theq} is small enough to have that the local diffeomorphism $\uu$ is defined in the ball of radius $2 \delta$ and center at the origin. Also, to simplify the notations we do not write $\bar U$, $\bar X^-$ and $\bar X^0$, but just $U$, $X^-$ and $X^0$. 
\subsubsection{A priori estimates}
\label{s:nont:zero}
We rewrite system \eqref{e:us:theq:b} as 
\begin{equation}
\label{e:us:spli}
\left\{
\begin{array}{ll}
            d X^- / d \tau = f^- (X^-, \, X^0) \\
            d X^0 / d \tau = f^0 (X^-, \, X^0) \\
\end{array}
\right.    
\end{equation}
The subspaces $\{ X^- = \vec 0 \}$ and $\{ X^0 = \vec 0 \}$ are locally invariant for \eqref{e:us:spli} since they represent respectively a center and the stable manifold. Thus, $f^- ( \vec 0, \, X^0 ) \equiv \vec 0$ for every $X^0$ and $f^0 ( X^-, \, \vec 0) \equiv 0$ for every $X^-$. As a consequence, 
\begin{equation}
\label{e:us:f:min}
          f^- (X^-, \, X^0) = A^- (X^-, \, X^0) X^-  \qquad  f^0 (X^-, \, X^0) = \hat{A}^0 (X^-, \, X^0) X^0
\end{equation}
for a suitable matrices $A^-$ and $A^0$. By construction, $A^-(\vec 0, \, \vec 0)$ admits only eigenvalues with strictly negative real part and $\hat{A}^0 (\vec 0, \, \vec 0)$ has only eigenvalues with zero real part. As a consequence, the following holds. Let $n_-$ denote the dimension of $X^-$ and fix a constant $c>0$ satisfying $Re \lambda < - c$ for every $\lambda$ eigenvalue of $A^- (\vec 0, \, \vec 0)$. Then there exists a constant $C_- >0$ such that 
\begin{equation}
\label{e:us:Cmeno}
          \forall \, \underline X^- \in \mathbb R^{n_-}, \qquad |e^{A^- (\vec 0, \, \vec 0) t   }  \underline X^- | \leq C_- e^{-c t} |   \underline X^-  |. 
\end{equation}
Also, if  $\delta$ is small enough and $|X^-(0)| < \delta$, then the solution of the Cauchy problem
$$
    \left\{
    \begin{array}{ll}
              d X^-/ d \tau = f^- (X^-, \, \vec 0) \\
              X^-(\tau = 0) = X^-(0)
    \end{array}
    \right.
$$
satisfies 
$$
    |X^- (\tau)| \leq C_- e^{- c \tau / 2 } |X^-(0)|,
$$
 where $c>0$ is as before a constant such that $Re \lambda < - c $ for every $\lambda$ eigenvalue of  $A^-(\vec 0, \, \vec 0)$.

Plugging \eqref{e:us:f:min} in \eqref{e:us:spli} we get
\begin{equation}
\label{e:us:theq:split}
          \left\{
          \begin{array}{ll}
            d X^- / d \tau = A^- (X^-, \, X^0) X^- \\
            d X^0 / d \tau = \hat{A}^0 (X^-, \, X^0) X^0. \\
\end{array}
\right.    
\end{equation}

In view of the applications discussed in Section \ref{s:invman} it is convenient to take into account the following situation.  Assume that there exists a continuously differentiable manifold $\mathcal Z_0$ containing the stable manifold $\{ X^0 = \vec 0 \}$ and satisfying 
\begin{equation}
\label{e:us:zero}
          f^0 (X^-, \, X^0) = \vec 0 \qquad \forall \, (X^-, \, X^0) \in \mathcal Z_0.
\end{equation}
Actually, this assumption is not restrictive, in the sense explained in Remark \ref{r:zero} at the end of Section \ref{s:nont:zero}. 

Applying, if needed, a local diffeomorphism, we can assume that $X^0 = (\zz, \, u_0)$ and that $\mathcal Z_0 = \{  \zz = \vec 0 \}$. Since the stable manifold is entirely contained in $\mathcal Z_0$, such a diffeomorphism does not produce any change on $X^-$, but only on $X^0$. 
In the following we assume that the constant $\delta$ in Hypothesis \ref{h:cutoff} is so small that the local diffeomorphism is defined in the ball of radius $2 \delta$ and center at the origin. 
 
Consider the system restricted on the center manifold $\{ X^ - =\vec 0 \}$: since the subspace $\{ \zz = \vec 0 \}$ is entirely made by equilibria, then we get that 
the equation
$$
     d X^0 / d \tau = \hat{A}^0 (\vec 0, \, X^0) X^0
$$
becomes 
\begin{equation}
\label{e:us:split2}
\left\{
\begin{array}{ll}
           d \zz / d \tau = \hat{B} (\vec 0, \, \zz, \, u_0) \zz  \\
           d u_0  / d \tau = \hat{C} (\vec 0, \, \zz, \, u_0 ) \zz, \\ 
\end{array}
\right.
\end{equation}
where $\hat{B}$ and $\hat{C}$ are suitable matrices. Note that, by construction, $\hat{B}(\vec 0, \, \vec 0, \, \vec 0)$ admits only eigenvalues with zero real part.  Fix a constant $\ee$ such that $Re \lambda < - \ee < 0$ for any $\lambda$ eigenvalue of $A^- (\vec 0, \, \vec 0)$: also, we impose $\ee < c$, where $c$ is the same as in \eqref{e:us:Cmeno}. Assuming that the constant $\delta$ in Hypothesis \ref{h:cutoff} is sufficiently small we can assume that every solution $\zz$ of \eqref{e:us:split2} satisfies 
\begin{equation}
\label{e:f:es:z}
    |\zz (\tau) | \leq \unpo e^{\ee |\tau| } | \zz (0) | 
\end{equation}
for some suitable constant $\unpo$.  Since in  \eqref{e:us:split2} the matrix $\hat{C}$ is uniformly bounded, we get that 
\begin{equation}
\label{e:f:es:u}
    |u_0 (\tau) - u_0 (0) | \leq \unpo e^{\ee | \tau|  } | \zz (0) | 
\end{equation}
for a constant $\unpo$ (possibly different from the one in \eqref{e:f:es:z}). We introduce the following notation: given a point $\underline X^0 = (\zz, \, u_0)$ on the center manifold we call $\underline Y^0$ the point 
\begin{equation}
\label{e:us:yo} 
         \underline Y^0 = (\vec 0, \, u_0). 
\end{equation}         
Clearly $\underline Y^0$ depends on $\underline X^0$, but to simplify the notations we won't express this dependence explicitly.
Combining \eqref{e:f:es:z} and \eqref{e:f:es:u} we then obtain 
\begin{equation}
\label{e:f:es:Y}
             |X^0 (\tau) - \underline Y^0 (0) | \leq k_0 e^{\ee |\tau|} | \zz (0) |      
\end{equation}
for a suitable constant $k_0$. 

Finally, note that, since both $A^-$ and $\hat{A}$ are zero when $|(X^-, \, X^0)| \ge 2 \delta$, then any non constant solution of \eqref{e:us:theq:split}
satisfies 
\begin{equation}
\label{e:glo:bd:0}
           | X^0 (\tau)  |  \leq 2 \delta  \qquad |X^-(\tau) | \leq 2 \delta \qquad  \forall \, \tau . 
\end{equation}
\begin{rem}
\label{r:zero}
The hypothesis that there exists a manifold of zeroes $\mathcal Z_0$ is not restrictive. Indeed, assume that the set of the zeroes of $f^0$ coincides with the stable manifold $\{ X^0 = \vec 0\}$. In this case, we can set $\zz = X^0$, the component $u_0$ disappears and given $\underline X^0$ the element $\underline Y^0$ is just $\underline X^0$ itself. This notation ensures that the estimate \eqref{e:f:es:Y} still holds. As it will be clear from the the following, the only fact about $\mathcal Z_0$ we exploit in the proof of Proposition \ref{p:exp:d} is estimate \eqref{e:f:es:Y}. As a consequence, Proposition  \ref{p:exp:d} can be extended to the case $\mathcal Z_0$ is just the stable manifold.

In other words, the presence of a manifold of zeroes wider then the stable manifold is not strictly necessary for the existence of the uniformly stable manifold introduced in Theorem \ref{p:unif}. However, it allows to get a sharper estimate in \eqref{e:p:exp:dec2}.
\end{rem}
\subsubsection{Linear change of variables}
\label{s:us:sec:ch}
In the statement of the following lemma we denote by $n_c$ the dimension of $X^0$, then $N=n_c + n_-$.  The proof is standard, so we omit it. 
\begin{lem}
\label{l:lin}
            For every $M > 0$, there exists a linear change of variables $\mathbb R^{n_c} \to \mathbb R^{n_c}$ such that in the new coordinates
            $X^0$ satisfies
            \begin{equation}
            \label{e:us:f:zero}
                      d X^0 / d \tau = \hat{A} (X^-, \, X^0) X^0,
            \end{equation}
            for a suitable matrix such that 
             \begin{equation}
            \label{e:hat0}
                \hat A^{0} (\vec 0, \, \vec 0) = \bar A^{0}  + N^0,
          \end{equation}
           where $ \bar A^0$ and $N^0$ enjoy the following properties:
            \begin{equation}
            \label{e:exp:dec:c}
                 |e^{\bar A^{0}t  } X^0 | \leq |X^0| \quad \forall t >0, \; X^0 \in \mathbb R^{n_c}
            \end{equation}
            and 
            \begin{equation}
            \label{e:nihil}
                   |N^{0} X^0 | \leq \frac{1}{M} |X^0|, \; X^0 \in \mathbb R^{n_c}.  
            \end{equation}
\end{lem}
We specify in the following how we chose the constant $M$.
\begin{rem}
\label{r:doppio}
           If we apply the linear change of coordinates introduced in Lemma \ref{l:lin}, then it is no more true that $X^0 =( \zz, \, u_0)$ where $\{ \zz = \vec 0 \}$ is the manifold $\mathcal Z_0$ of equilbria for $f^0$. However, estimate \eqref{e:f:es:Y} still holds, provided that we change if needed the value of the constant $k_0$ and we take, instead of $X^0(\tau)$, $\underline Y_0$  and $\zz (0)$, their images trough the linear change of variables.  
\end{rem}

\subsection{Uniformly stable manifold of an orbit}
\label{s:us:orbit}
We are now ready to introduce Theorem \ref{p:unif}. In formula \eqref{e:p:exp:dec2}, $\zz$ is the component of $X^0 =( \zz, \, u_0)$ according to the decomposition introduced in Section \ref{s:us:first:ch}.
\begin{teo}
\label{p:unif}
            Let Hypotheses \ref{h:cutoff} and \ref{h:neg} hold. If the constant $\delta$ in Hypothesis \ref{h:cutoff} is sufficiently small, then the following holds.

            Fix an orbit $Y^0 (\tau) = \big(\vec 0, \, X^0(\tau) \big)$ of 
            \begin{equation}
            \label{e:teo:sy}
          \left\{
          \begin{array}{ll}
            d X^- / d \tau = A^- (X^-, \, X^0) X^- \\
            d X^0 / d \tau = \hat{A}^0 (X^-, \, X^0) X^0 \\
\end{array}
\right.    
\end{equation}
           that lies on the center manifold and satisfies $|X^0(0)| < \delta$. Then we can define a uniformly stable manifold relative to $Y^0(\tau)$. This manifold 
            is defined in the ball of radius $ \delta$ and center at the origin, is parameterized by  $\{ X^0 = \vec 0\}$ and is tangent to this subspace at the origin. Also, it is locally invariant   for 
            \eqref{e:teo:sy}, meaning that if the initial datum lies on the manifold, then $\big( X^-(\tau), \, X^0(\tau) \big)$ belongs to the uniformly stable manifold for $|\tau |$ sufficiently small.  
            Every orbit lying on the uniformly stable manifold relative to $Y^0(\tau)$ can be decomposed as 
            \begin{equation}
            \label{e:decompo}
                 X(\tau) = Y^0 (\tau) + Y^{-} (\tau) + U^p (\tau),
            \end{equation}
            where the components  $Y^{-}= \big(X^-(\tau), \vec 0 \big)$ and $U^p(\tau)$ satisfy respectively 
            \begin{equation}
            \label{e:p:exp:decst}
                           |X^- (\tau) | \leq k_-  |X^- (0)| e^{-c \tau/ 2}
            \end{equation}
            and 
            \begin{equation}
            \label{e:p:exp:dec2}
                 |U^p (\tau) | \leq k_p |\zz (0)| \, |X^- (0)| e^{-c \tau/ 4}.
            \end{equation}
            In \eqref{e:p:exp:decst} and \eqref{e:p:exp:dec2},  $c$, $k_-$ and $k_p$ are suitable constants. 
            In particular, $c$ is the same as in \eqref{e:us:Cmeno}. 
            \end{teo} 
\subsection{Proof of Theorem \ref{p:unif}}
\label{s:us:proof}
Let the orbit $(\vec 0, \, X^0(\tau))$ be given. We denote by $
           \underline X^0 = X^0(0) $
and by $\underline Y^0$ the corresponding projection, defined as in \eqref{e:us:yo}. Also, if we write $\underline X^0 = (\zz(0), \, u_0(0))$ then we set
\begin{equation}
\label{e:us:z0}
           \underline \zz = \zz (0).
\end{equation}
By definition,  $X^0(\tau)$ is a solution of the Cauchy problem
\begin{equation}
\label{e:us:cau:x0}
\left\{
\begin{array}{ll}
           d X^0 / d \tau = \hat A^0 (\vec 0, \, X^0) \\
           X^0 (0) = \underline X^0 \\
\end{array}
\right.
\end{equation}
The proof of Theorem  \ref{p:unif}  is divided in several steps: in Section \ref{s:us:fsp} we introduce the spaces of functions we exploit in the proof. In Section \ref{s:us:stable} we are concerned with the component $Y^-(\tau)$ in \eqref{e:decompo}. In Section \ref{s:us:pert} we deal with the ``perturbation" term $U^p(\tau)$ in \eqref{e:decompo}. Both the components $Y^-(\tau)$ and $U^p(\tau)$ are obtained as fixed points of suitable contractions: in Section \ref{s:us:fre} we study their regularity applying Lemma \ref{l:fre}. Finally, in Section \ref{e:us:con} we conclude the proof of Theorem \ref{p:unif} putting together all the considerations carried on in Sections \ref{s:us:fsp}. \ref{s:us:stable}, \ref{s:us:pert} and \ref{s:us:fre}.   
\subsubsection{Definition of the functional spaces}
\label{s:us:fsp}
Let $n_-$ denote the dimension of $X^-$. Also, $n_c$ denotes the dimension of $X^0$, as in the statement of Lemma \ref{l:lin}. 

In the following we exploit the following Banach spaces of functions:
\begin{equation}
\label{e:us:sp:y-}
    \mathcal Y^- = \Big\{ 
                       X^- \in \mathcal C^0 \big( [0, \, + \infty [ , \; \mathbb R^{n_-}  \;  \big): \; \|  X^-  \|_-    <   + \infty 
                       \Big\} 
\end{equation}
and 
\begin{equation}
\label{e:us:sp:yo}
    \mathcal Y^0 = \Big\{ 
                       X^0 \in \mathcal C^0 \big( [0 , \, + \infty [ , \; \mathbb R^{n_c}  \;  \big): \; \|  X^0  \|_0   <  + \infty   
                       \Big\}, 
\end{equation}
where the norms $\| \cdot \|_{- }$ and $\|  \cdot \|_{0}$ are defined as follows:
\begin{equation}
\label{e:us:sp:norms}
      \|  X^-  \|_-  = \sup_{\tau} \big\{ e^{ c \tau / 2 } | X^-(\tau) | \big\} \qquad 
         \|  X^0  \|_0  = \sup_{\tau} \big\{ e^{ - \ee |\tau | } | X^0(\tau) |  \big\} 
\end{equation}         
The constants $c$ and $\ee$ are the same as in \eqref{e:us:Cmeno} and  \eqref{e:f:es:Y} respectively. Also, we consider the following closed subsets of $Y^-$ and $Y^0$:
\begin{equation}
\label{e:st:y:d}
    \mathcal Y^-_{\delta} = \Big\{ 
                       X^- \in \mathcal C^0 \big( [    0 , \, + \infty [ , \; \mathbb R^{n_-}  \;  \big): \; \|  X^-  \|_-    \leq    k_- \delta.
                       \Big\} \qquad  \mathcal Y^0_{\delta} = \Big\{ 
                       X^0 \in \mathcal C^0 \big( [   0 , \, + \infty [ , \; \mathbb R^{n_c}  \;  \big): \; \|  X^0  \|_0    \leq  k_0\delta.
                       \Big\} 
\end{equation}
We specify in the following how to determine the exact value of the constant $k_-$, while the constant $k_0$ is the same as in \eqref{e:f:es:Y}. Also, note that the spaces $\mathcal  Y^-_{\delta}$ and $\mathcal Y^0_{\delta}$ are equipped with the same norms $\| \cdot \|_-$ and $\| \cdot \|_0$ as $\mathcal  Y^-$ and $\mathcal Y^0$ respectively.  

We will also need the space of functions defined as follows. Let $c$ be as \eqref{e:us:Cmeno} and let $a \in [0, \, c[$. Consider the space
\begin{equation}
\label{e:us:sp:yp}
    \mathcal Y^p_a = \Big\{ 
                       (U^-, \, U^0) \in \mathcal C^0 \big( [  0 , \, + \infty [ , \; \mathbb R^{n_c+ n_-}  \;  \big): \; \|  X^-  \|_{pert}   <   + \infty 
                       \Big\} 
\end{equation}
which depends on $a$ because it is equipped with the norm 
$$
     \| (U^-, \, U^0)  \|_{pert}  =  \sup_{\tau} \Big\{ e^{(c+a) \tau / 4 }\Big[   |( U^- (\tau)|+ |U^0 (\tau) |)  \Big] \Big\}.
$$
Also, we will exploit the closed subset  
\begin{equation}
\label{e:us:sp:ypd}
    \mathcal Y^p_{\delta a} = \Big\{ 
                       (U^-, \, U^0) \in \mathcal C^0 \big( [    0 , \, + \infty [ , \; \mathbb R^{n_c+ n_-}  \;  \big): \; \|   (U^-, \, U^0) \|_{pert}    \leq   k_p \delta^2
                       \Big\}, 
\end{equation}
which is equipped with the same norm as $\mathcal Y^p$. We specify in the following how to determine the values of the constants $k_p$ and $a$.  
\subsubsection{Analysis of the stable component}
\label{s:us:stable}
This step is devoted to the definition of $Y^{-} (\tau) = ( X^- (\tau), \, \vec 0 )$. Fix a vector $\underline X^- \in  \mathbb R^{n_-}$ satisfying $|\underline X^- | < \delta$. 

We define $X^-(\tau)$ as the solution of the Cauchy problem 
\begin{equation}
\label{e:+:stac}
\left\{
\begin{array}{ll}
     d X^- / d \tau = A^{-} (X^-, \, \underline Y^0 ) X^- \\
     X^- (0) = \underline X^-,
\end{array}
\right.
\end{equation}
where $\underline Y^0$ is given by \eqref{e:us:yo}. 
It is known that, for any fixed  $\underline Y^0$ and $\underline X^-$,  $X^-$ can be obtained as the fixed point of the application
$$
     T^-:  \mathcal Y^-_{\delta} \to \mathcal  Y^-_{\delta}
$$
defined by
\begin{equation}
\label{e:st:con}
           T^- (X^-) [\tau] =  e^{\bar A^- \tau} \underline X^-  + \int_0^{\tau} e^{\bar A^- (\tau-s)} \Big[  A^- \big(X^-(s), \, \underline Y^0 \big) - \bar A^-    \Big] X^- (s) ds 
\end{equation}
where $\bar A^- = A^- (\vec 0, \, \vec 0)$. The space $\mathcal  Y^-_{\delta}$ is defined in \eqref{e:st:y:d}. More precisely, if the constant $k_-$ in \eqref{e:st:y:d} satisfies $k_- \leq C_-$ and the constant $\delta$ is \eqref{e:st:y:d} is sufficiently small, then the map $T^-$ takes values in $\mathcal  Y^-_{\delta}$ and is indeed a contraction. Also, the fixed point satisfies
\begin{equation}
\label{e:p:bd-}
     |X^- (\tau) | \leq k_- |\underline X^-| e^{-c \tau /2 } 
\end{equation}
We are now interested in the differentiability of the fixed point with respect to $\underline Y^0$ and $\underline X^-$. To study it, we recall that
$$
     \mathcal Z_0 = \{ (\underline X^-, \, \vec 0, \, \underline u_0)  \} \subseteq \mathbb R^N 
$$
We then regard $T^-$ as an application 
$$
     T^- : \mathcal Z_0  \times \mathcal  Y^-_{\delta} \to \mathcal  Y^-
$$
and we verify that the hypotheses of Lemma \ref{l:fre} are satisfied.  The space $\mathcal  Y^-$ is defined by \eqref{e:us:sp:y-}. 
The Frech\'{e}t derivative of $T^-$ with respect to $(\underline X^-, \, \underline Y^0)$ is a linear map $\mathcal T^- \in \mathcal L (\mathcal Z_0, \, \mathcal  Y^-)$. Evaluated at the point $(\underline h^-, \, \underline h^0) \in \mathcal Z_0$ it takes the value 
$$
     \mathcal T^-     (h^-, \, h^0) [\tau] =  e^{\bar A^- \tau} \underline h^- +  \int_0^{\tau} e^{\bar A^- (\tau-s)} \Big[ D_{\underline Y^0} 
     A^- \big(X^-(s), \, \underline Y^0 \big)  [\underline h^0]\Big] X^- (s) ds 
$$
 In the previous expression, $D_{\underline Y^0} 
     A^- \big(X^-(s), \, \underline Y^0 \big)  [\underline h^0]$ denotes the differential of the function $A^- \big(X^-(s), \, \underline Y^0 \big) $ with respect $\underline Y^0$, applied to the vector $\underline h^0$. If $\underline X^- = \vec 0$ then, no matter what $\underline Y^0$ is, the differential  $\mathcal T^-$ maps $(\underline h^-, \, \underline h^0)$ into the function $ e^{\bar A^- \tau} \underline h^- $.   

The Frech\'{e}t derivative of $T^-$ with respect to $X^-$ is a linear map $\mathcal S^- \in \mathcal L (\mathcal  Y^-, \, \mathcal  Y^-)$. Evaluated at the point $ h^- \in \mathcal  Y^-$ it takes the value 
$$
     \mathcal S^-     (h^-) [\tau] =   \int_0^{\tau} e^{\bar A^- (\tau-s)}  \Big\{ 
     \Big[  A^- \big(X^-(s), \, \underline Y^0 \big) - \bar A^-  \Big] h^- (s) + 
     \Big[ D_{\underline X^-} 
     A^- \big(X^-(s), \, \underline Y^0 \big)  [h^- (s) ]\Big] X^- (s) \Big\} ds 
$$
 In the previous expression, $D_{X^-} 
     A^- \big(X^-(s), \, \underline Y^0 \big)  [ h^- (s)]$ denotes the differential of the function $A^- \big(X^-(s), \, \underline Y^0 \big) $ with respect $X^-$, applied to the vector $h^-(s)$. 
          
Both  $\mathcal T^-$ and $\mathcal S^-$ are continuous if viewed as maps from $\mathcal Z_0 \times \mathcal Y_{\delta}$ to $ \mathcal L (\mathcal Z_0, \, \mathcal Y^-)$ and  $\mathcal L (\mathcal Y^-, \, \mathcal Y^-)$ respectively. Thus, the hypotheses of Lemma \ref{l:fre} are verified and hence the application
$$
     \mathcal Z_0  \to \mathcal Y^-_{\delta}
$$       
which associates to $(\underline X^-, \, \underline Y^0)$ the fixed point of \eqref{e:st:con} is continuously differentiable (in the sense of Frech\'{e}t). When both $\underline X^- = \vec 0$ and $\underline Y^0 = \vec 0$ the Frech\'{e}t derivative is the functional that maps $(\underline h^-, \, \underline h^0) \in \mathcal Z_0$ into the function $ e^{\bar A^- \tau} \underline h^- $. 
\subsubsection{Analysis of the component of perturbation}
\label{s:us:pert}
This step is devoted to the definition the component $U^p (\tau)$. First, we apply the change of variables introduced in Lemma \ref{l:lin} and we get that the matrix $\hat{A} (X^-, \, X^0)$ in \eqref{e:us:theq} satisfies 
$$
    \hat{A} (\vec 0, \, \vec 0) = \bar A^0 + N^0,
$$
where $\bar A^0$ and $N^0$ enjoy \eqref{e:exp:dec:c} and \eqref{e:nihil} respectively. Relying on Remark \ref{r:doppio}, we can still exploit estimate \eqref{e:f:es:Y}.

We impose that $ X(\tau) = Y^0 (\tau) + Y^{st} (\tau) + U^p (\tau)$ is a solution of \eqref{e:us:theq:split}. We then write $U^p (\tau) = \big( U^-, \, U^0)^T$ and, subtracting \eqref{e:us:cau:x0} and \eqref{e:+:stac} from \eqref{e:us:theq:split}, we get 
\begin{equation}
\label{e:pert}
\left\{
\begin{array}{llll}
            \displaystyle{ d U^- / d \tau = \bar A^-  U^- + \Big[ A^- (X^- + U^-, \, X^0 + U^0 ) -  A^- (\vec 0, \, \vec 0)  \Big] U^- } \\
            \qquad + 
            \displaystyle{ \Big[ A^- (X^- + U^-, \, X^0 + U^0 ) -  A^- (X^-, \, \underline Y^0)  \Big] X^-  } \\ 
             \displaystyle{ d U^0 / d \tau = \bar A^0  U^0 + N^0 U^0 + \Big[ \hat{A}^0 (X^- + U^-, \, X^0 + U^0 ) -  \hat A^0 (\vec 0, \, \vec 0)  \Big] U^0 }\\
           \qquad   + 
           \displaystyle{\Big[ \hat{A}^0(X^- + U^-, \, X^0 + U^0 ) -  \hat{A}^0 ( \vec 0, \, X^0)  \Big] X^0 } \\ 
\end{array}
\right.
\end{equation}
Here, $\bar A^- =  A^-(\vec 0, \vec 0)$. 

Let $\mathcal Y^p_{\delta a}$ be the metric space \eqref{e:us:sp:ypd} and consider the application $T_p$, defined for $(U^-, \, U^0) \in \mathcal Y^p_{\delta a}$
 as follows:
\begin{equation}
\label{e:p:con:def}
\begin{split}
&             T^1_p (U^-, \, U^0 ) [\tau]= \displaystyle{  \int_0^{\tau} e^{\bar A^- (\tau-s)} \Big\{ \big[ A^-\big( X^- (s)+ U^-(s), \, X^0 (s)+ U^0 (s) \big) - A^- \big( X^-(s), \, \underline Y^0 \big)    \big]      X^- (s) } \\
&          \qquad \qquad \qquad \qquad +
             \big[ A^-\big( X^- (s)+ U^-(s), \, X^0 (s)+ U^0 (s) \big) - A^- \big( \vec 0, \, \vec 0 \big)    \big] U^- (s) \Big\} ds   \\
&         T^2_p (U^-, \, U^0 ) [ \tau] = \int_{+ \infty }^{\tau} e^{ \bar A^0 (\tau -s ) }
             \Big\{  \Big[ N^0 + \hat A^0 \big(X^- (s) + U^- (s), \, X^0 (s)+ U^0  (s) \big)  - \hat A^0 \big( \vec 0, \, \vec 0 \big) \Big] U^0 (s)  \\
&         \qquad    \qquad \qquad \qquad +
             \Big[ \hat A^0 \Big(X^- (s) + U^- (s), \, X^0 (s)+ U^0  (s) \big)  - \hat A^0 (\vec 0, \, X^0 (s) \big) \Big] X^0 (s)     \Big\} ds \\
\end{split}
\end{equation}
In the previous expression, $X^-$ is the solution of \eqref{e:+:stac} and $X^0$ is the solution of \eqref{e:us:cau:x0}. We want to show that $T_p$ maps $\mathcal Y^p_{\delta a}$ into itself, provided that $\delta$ is sufficiently small. 
We have 
\begin{equation}
\label{e:est:p1}
\begin{split}
&          |T^1_p (U^-, \, U^0 ) [\tau] | 
         \leq \int_0^{\tau}C_-  e^{- c(  \tau -s )  } \Big\{ L \Big[   |U^- (s)| + |U^0(s) | + | X^0 (s)  - \underline Y^0|   \Big]  |X^- (s) |  \\
&         \qquad  +
           L \Big[   |X^- (s)| + |U^-(s) | + | X^0 (s)  |   + |U^0 (s) | \Big]  |U^- (s) |   \Big\} ds  \\
&         \leq C_- e^{- c \tau} \int_0^{\tau} e^{c s }  L \Big[ 2   k_p \delta^2 e^{-s (c + a)/4}  + k_0 |\underline \zz | e^{\ee s} \Big] 
            k_- |\underline X^-|  e^{-cs/ 2}  \\
 &        \qquad
             +   C_- e^{- c \tau } \int_0^{\tau} e^{c s   }   L \Big[   k_- |\underline X^-|  e^{-cs/2}   + 2 k_p   \delta^2     e^{-s(c +a )/4}  + 2 \delta   \Big] 
            k_p  \delta^2   e^{-s(c+a)/ 4}  ds \\
&       \leq   \Big[  \frac{8}{c - a } C_- L k_p k_-  \delta  \Big] \delta^2 e^{-  \tau (3c + a )/4 } +  \frac{2}{c + 2 \ee }  L C_-  k_0 k_-  \delta^2 e^{-\tau ( 2c - 4 \ee) /4 } +  \Big[   \frac{4}{c - a }  
         L  C_-  k_- k_p \delta \Big]\delta^2  e^{-  \tau (3c + a )/4 }  \\
 &    \quad         + 
          \Big[  \frac{4}{c- a }  L C_-  k_p  \delta^2  \Big]  k_p   \delta^2    e^{-  \tau (c + a )/ 2 }    +       \Big[  \frac{8}{3c - a } \delta  \Big] L C_- k_p   \delta^2    e^{-  \tau (c + a) / 4 }.        
\end{split}
\end{equation}
In the previous expression, $C_-$ is the same constant as in \eqref{e:us:Cmeno} and $L$ is a Lipschitz constant of $A^-(X^-, \, X^0)$ with respect to both the variables $X^-$ and $X^0$. To obtain \eqref{e:est:p1} we exploit \eqref{e:f:es:Y}, \eqref{e:glo:bd:0},  \eqref{e:p:bd-} and the fact that, belonging to $\mathcal Y^p_{\delta a}$, $(U^-, \, U^0)$ satisfies 
\begin{equation}
\label{e:est:u}
    |U^- (\tau)|, \; |U^0 (\tau)| \leq k_p \delta^2 e^{- \tau ( c + a) / 4}. 
\end{equation}
Also, the term $\underline \zz$ is the same as in \eqref{e:us:z0} and we rely on the fact that $|\underline X^-|$, $|\underline \zz| < \delta$.

In the following expression $L$ denotes a Lipschitz constant of $\hat{A}^0(X^-, \, X^0)$ with respect to both the variables $X^-$ and $X^0$. Also, we exploit the estimates \eqref{e:exp:dec:c}, \eqref{e:nihil}, \eqref{e:p:bd-}, \eqref{e:est:u} and \eqref{e:glo:bd:0}.  
\begin{equation}
\label{e:est:p2}
\begin{split}
&        | T^2_p (U^-, \, U^0 ) [ \tau]  | \leq \int_{+ \infty}^{\tau} |N^0 U^0 (s) | +L \Big[   |X^-(s) | + |U^- (s) | + |X^0 (s) | + |U^0 (s) |  \Big] |U^0 (s)| ds \\
&       \quad +   \int_{+ \infty}^{\tau} L  \Big[   |X^-(s) | + |U^- (s) | +  |U^0 (s) |  \Big] |X^0 (s)| ds \\ 
&       \leq  \int_{+ \infty}^{\tau} \frac{1}{M}  |U^0 (s)| ds  + L  \int_{+ \infty}^{\tau}  
           \Big[ k_- |\underline X^-| e^{-c s /2} +2  k_p \delta^2  e^{- s (c+a)  /4} + 2 \delta  \Big]  k_p \delta^2  e^{-s (c +a ) /4} ds  \\
&       \quad  +  L   k_-    \int_{+ \infty}^{\tau}   |\underline X^- |   e^{-cs / 2 }2 \delta ds + L    \int_{+ \infty}^{\tau} 2   k_p \delta^2 e^{- s (c +a) /4}  2  \delta ds     \\
&          \leq \frac{4}{M (c + a) } k_p \delta^2   e^{- \tau (c +a)/ 4} + \frac{4 L  k_- \delta  }{3c + a }  k_p \delta^2  e^{-\tau (3 c +a ) /4}  + 
            \frac{4 L k_p \delta^2}{c + a } k_p \delta^2   e^{-\tau (c +a ) /2 } + 
            \frac{8 L  \delta}{c + a } k_p \delta^2   e^{-\tau (c +a ) / 4 } \\ 
&        \quad + \frac{4 }{c} k_- L \delta^2  e^{- c \tau / 2 } +
           \frac{16 L \delta }{c + a } k_p \delta^2 e^{- \tau (c +a) /4}  \\
\end{split}
\end{equation}
Combining \eqref{e:est:p1} and \eqref{e:est:p2} we get the following. Assume that the constant $k_p$ in \eqref{e:us:sp:ypd} is sufficiently large (namely, $k_p \ge 4 L k_- /c$). Then for every $a \leq c - 4 \ee$ we can choose $\delta$ and $M$ in such a way that $T_p$ take values into $\mathcal Y^p_{\delta a}$.  Also, estimates similar to \eqref{e:est:p1} and \eqref{e:est:p2} ensure that one can choose the constants in such a way that $T_p$ is a strict contraction. As a remark,  we point out that, the bigger is $a$, the smaller is $\delta$.  

We set $a =12 \ee$ and we choose $\delta$ in such a way that $T^p$ is a contraction from $\mathcal Y^p_{\delta 12 \ee}$ to itself. The constant $\ee$ is the same as in \eqref{e:f:es:z}. 
However, in the following we regard $T^p$ as a map $\mathcal Y^p_{\delta 0} \to\mathcal Y^p_{\delta 0}$, where $\mathcal Y^p_{\delta 0}$ is the space \eqref{e:us:sp:ypd} obtained setting $a=0$. In this way, we obtain that  $T^p$ is a contraction on $\mathcal Y^p_{\delta 0}$, but, thanks to our choice of $\delta$, the fixed point automatically satisfies the sharper estimate 
\begin{equation}
\label{e:imp:fp} 
          |U^- (\tau) |, \; | U^0 (\tau) | \leq k_p \delta^2 e^{ - \tau ( c + 12 \ee)/ 4}.  
\end{equation}
Also, in the definition of the space $\mathcal Y^p_{\delta 0}$ one can take $\delta^2 = |\underline \zz| \, |\underline X^-|$ and hence 
\begin{equation}
\label{e:us:es:u}
      |U^- (\tau) |, \; | U^0 (\tau) | \leq k_p e^{ - \tau ( c + 12 \ee)/ 4}  |\underline \zz| \, |\underline X^-|,
\end{equation}
where $\underline X^-$ is defined by  \eqref{e:+:stac}. Also, to simplify the notations in the previous expression we denote by $\underline \zz$ the point obtained applying the change of coordinates introduced in Lemma \ref{l:lin} to the vector $( \underline \zz, \, \vec 0)$ defined by \eqref{e:us:z0}. 
 \subsubsection{Frech\'{e}t differentiability of the component of perturbation}
 \label{s:us:fre} 
We are now concerned with the Frech\'{e}t differentiability of the fixed point of the map $T_p$ defined by \eqref{e:p:con:def}. Since $T_p (U^-, \, U^0)$ depends on $X^-$ and $X^0$, we regard $T_p$ as a map
\begin{equation}
\label{e:T}
     T:  \mathcal Y^p_{\delta 0}  \times  \mathcal Y \to  \mathcal Y^p_{\delta 0}.
\end{equation}
In the previous expression, $\mathcal Y=\mathcal Y^- \times \mathcal Y^0$, where $\mathcal Y^-$ and $\mathcal Y^0$ are defined by \eqref{e:us:sp:y-} and \eqref{e:us:sp:yo} respectively. Also, they satisfy $X^- \in \mathcal Y^-$ and $X^0 \in \mathcal Y^0$.

The proof of the differentiability relies on Lemma \ref{l:fre} (taking $Y = \mathcal Y$ and $X =  \mathcal Y^p_{\delta 0}$). We thus verify that the hypotheses of Lemma \ref{l:fre} are satisfied. 

To simplify the exposition, we write \eqref{e:p:con:def} as  
\begin{equation}
\label{e:p:con:ch}
\begin{split}
&             T^1_p (U^-, \, U^0 ) [\tau]= \displaystyle{  \int_0^{\tau} e^{\bar A^- (\tau-s)} \Big\{ F \Big( X^- (s), \; U^- (s), \, X^0 (s), \;  U^0 (s) \Big) X^- (s) } \\
&            \qquad \qquad \qquad 
                 + G  \Big( X^- (s) +  U^- (s), \, X^0 (s) + U^0 (s) \Big)  U^- (s)  \Big\} ds   \\
&         T^2_p (U^-, \, U^0 ) [ \tau] = \int_{+ \infty }^{\tau} e^{ \bar A^0 (\tau -s ) }
             \Big\{  \Big[ N^0 + H \Big( X^- (s) + U^- (s), \, X^0 (s) + U^0 (s) \Big)   \Big] U^0 (s)  \\
&         \qquad    \qquad \qquad +
             L \Big( X^- (s) + U^- (s), \, X^0 (s), \,  U^0 (s) \Big)    X^0 (s)     \Big\} ds, \\
\end{split}
\end{equation}
where the functions $F$, $G$, $H$ and $L$ satisfy
\begin{equation}
\label{e:p:bd:FGHL}
      F \Big( X^- (s), \; \vec 0, \, X^0 (s), \;  \vec 0  \Big)  \equiv 0 \quad  G  \Big(  \vec 0, \,  \vec 0 \Big) \equiv \vec 0
     \quad 
     H  \Big( \vec 0, \, \vec 0 \Big) \equiv \vec 0 \quad L \Big( \vec 0, \,  X^0, \,  \vec 0 \Big) \equiv \vec 0.
\end{equation}
Note that $X^0(s) \equiv \underline Y^0$ is an equilibrium for \eqref{e:us:cau:x0}.  

Relying on \eqref{e:imp:fp}, one can show that the condition \eqref{e:us:suff:Lip} is verified here, so applying Remark \ref{r:lip} we get that the fixed point $(U^-, \, U^0)$ is Lipschitz continuous with respect to $(X^0, \, X^-)$.
 
Concerning the Frech\'{e}t  differentiability of $T_p$ with respect to $(X^0, \, X^-)$, we proceed as follows. Fix an element ${(U^0, \, U^-, \, X^0, \, X^-) \in \mathcal Y^p_{\delta 0} \times \mathcal Y}$ satisfying the estimates \eqref{e:f:es:Y}, \eqref{e:glo:bd:0}, \eqref{e:p:bd-}  and \eqref{e:imp:fp}.
 The Frech\'{e}t  differential of $T_p$ with respect to $(X^-, \, X^0)$ computed at the point $(U^0, \, U^-, \, X^0, \, X^-)$  is a linear map $\mathcal T 
\in \mathcal L (\mathcal Y, \,  \mathcal Y^p_{0})$. The image of the element $(h^-, \, h^0) \in \mathcal Y = \mathcal Y^- \times \mathcal Y^0$ is given by 
\begin{equation}
\label{e:p:diff:x}
\begin{split}
&             \mathcal T^1_p (h^-, \, h^0 ) [\tau]= \displaystyle{  \int_0^{\tau} e^{\bar A^- (\tau -s)} \Big\{ F \Big( X^- (s), \; U^- (s), \, X^0 (s), \; U^0 (s) \Big) h^- (s)} \\
&           \qquad \qquad \qquad 
                + \Big[ D_{X^-} F   \Big( X^- (s), \; U^- (s), \, X^0 (s), \; U^0 (s) \Big) h^- (s)  \Big] X^- (s) \\
&           \qquad \qquad \qquad                 +  
                \Big[ D_{X^-} G   \Big( X^- (s) + U^- (s), \, X^0 (s) + U^0 (s) \Big) h^- (s)  \Big] U^- (s)  \\
&                \qquad \qquad \qquad 
              + \Big[ D_{X^0} F   \Big( X^- (s), \; U^- (s), \, X^0 (s), \; U^0 (s) \Big) h^0(s)  \Big] X^- (s) \\
              &           \qquad \qquad \qquad +  
                \Big[ D_{X^0} G   \Big( X^- (s) + U^- (s), \, X^0 (s) + U^0 (s) \Big) h^0 (s)  \Big] U^- (s)   \Big\} ds   \\
&         \mathcal T^2_p (h^-, \, h^0 ) [ \tau] = \int_{+ \infty }^{\tau} e^{ \bar A^0 (\tau -s ) }
             \Big\{    \Big[ D_{X^-} H   \Big( X^- (s) +  U^- (s), \, X^0 (s) + U^0 (s) \Big) h^- (s)  \Big] U^0 (s) \\
&           \qquad \qquad \qquad                 +  
                \Big[ D_{X^-} L   \Big( X^- (s) +  U^- (s), \, X^0 (s), \;  U^0 (s) \Big) h^- (s)  \Big] X^0 (s)  \\
&                \qquad \qquad \qquad 
              + \Big[ D_{X^0} H   \Big( X^- (s) + U^- (s), \, X^0 (s) + U^0 (s) \Big) h^0(s)  \Big] U^0 (s) \\
 &           \qquad \qquad \qquad +  
              L  \Big( X^- (s) + U^- (s), \, X^0 (s), \; U^0 (s) \Big) h^0 (s) \\               
              &           \qquad \qquad \qquad +  
                \Big[ D_{X^0} L   \Big( X^- (s) + U^- (s), \, X^0 (s), \, U^0 (s) \Big) h^0 (s)  \Big] X^0(s)   \Big\} ds   \\
\end{split}                
\end{equation}
In the previous expression, $\big[ D_{X^-} F   \big( X^- (s), \; U^- (s), \, X^0 (s), \; U^0 (s) \big) h^- (s)  \big] $ denotes the differential of the matrix valued function $F$ with respect to the variable $X^-$. The differential is computed at the point $ \big( X^- (s), \; U^- (s), \, X^0 (s), \;  U^0 (s) \big)$ and is applied to the vector $h^-(s)$.  
To prove that indeed
$$
   {\big( \mathcal T^1(h^-, \, h^0 ), \,  \mathcal T^2 (h^-, \, h^0 ) \big) \in \mathcal Y^p_0}
$$ 
one exploits estimate \eqref{e:imp:fp} and the identity $L( \vec 0, \, X^0, \, \vec 0) \equiv \vec 0$.

We now discuss the the Frech\'{e}t differentiability of $T_p$ with respect to $(U^0, \, U^-)$.  Fix an element ${(U^0, \, U^-, \, X^0, \, X^-) \in \mathcal Y^p_{\delta 0} \times \mathcal Y}$. The Frech\'{e}t  differential of $T_p$  with respect to $(U^0, \, U^-)$, evaluated at the point $(U^0, \, U^-, \, X^0, \, X^-)$, is a linear map $\mathcal S
\in \mathcal L ( \mathcal Y^p_0, \, \mathcal  Y^p_0)$ and the image of the element $(h^-, \, h^0) \in \mathcal Y^p_0$ is given by 
\begin{equation}
\begin{split}
&             \mathcal S^1(h^-, \, h^0 ) [\tau]= \displaystyle{  \int_0^{\tau} e^{\bar A^- (\tau -s)} \Big\{ 
                 \Big[ D_{U^-} F   \Big( X^- (s), \; U^- (s), \, X^0 (s), \;  U^0 (s) \Big) h^- (s)  \Big] X^- (s) } \\
&           \qquad \qquad \qquad  + G   \Big( X^- (s) +  U^- (s), \, X^0 (s) + U^0 (s) \Big) h^- (s)      \\
&           \qquad \qquad \qquad           +  
                \Big[ D_{U^-} G   \Big( X^- (s) +  U^- (s), \, X^0 (s) + U^0 (s) \Big) h^- (s)  \Big] U^- (s)  \\
&                \qquad \qquad \qquad 
              + \Big[ D_{U^0} F   \Big( X^- (s), \; U^- (s), \, X^0 (s), \,  U^0 (s) \Big) h^0(s)  \Big] X^- (s) \\
              &           \qquad \qquad \qquad +  
                \Big[ D_{U^0} G   \Big( X^- (s), \; U^- (s), \, X^0 (s) + U^0 (s) \Big) h^0 (s)  \Big] U^- (s)   \Big\} ds   \\
&         \mathcal S^2 (U^-, \, U^0 ) [ \tau] = \int_{+ \infty }^{\tau} e^{ \bar A^0 (\tau -s ) }
             \Big\{     \Big[ D_{U^-} H   \Big( X^- (s) + U^- (s), \, X^0 (s) + U^0 (s) \Big) h^- (s)  \Big] U^0 (s) \\
&           \qquad \qquad \qquad                 +  
                \Big[ D_{U^-} L   \Big( X^- (s) +  U^- (s), \, X^0 (s), \,  U^0 (s) \Big) h^- (s)  \Big] X^0 (s)  \\
&                \qquad \qquad \qquad + 
                N^0 h^0 (s) +  H   \Big( X^- (s), \; U^- (s), \, X^0 (s) + U^0 (s) \Big) h^0 (s)  \\
&          \qquad \qquad \qquad 
              + \Big[ D_{U^0} H   \Big( X^- (s) +  U^- (s), \, X^0 (s) + U^0 (s) \Big) h^0(s)  \Big] U^0 (s) \\            
              &           \qquad \qquad \qquad +  
                \Big[ D_{U^0} L   \Big( X^- (s) + U^- (s), \, X^0 (s), \, U^0 (s) \Big) h^0 (s)  \Big] X^0(s)   \Big\} ds   \\
\end{split}                
\end{equation}
One can verify that, if ${(U^0, \, U^-, \, X^0, \, X^-) \in \mathcal Y^p_{\delta 0}  \times \mathcal Y}$, then indeed $\mathcal S(h^-, \, h^0) \in \mathcal Y^p_0$. Also, $\mathcal S$ is continuous as  a map from $ X^p \times \mathcal Y$ in $\mathcal L ( \mathcal Y^p_0, \, \mathcal  Y^p_0)$. 

This shows that the hypotheses of Lemma \ref{l:fre} are all verified.

\subsubsection{Conclusion}
\label{e:us:con}
Applying Lemma \ref{l:fre}, we get that the map 
\begin{equation}
\label{e:p:fp:reg}
         \mathcal Y \to \mathcal Y^p_{\delta 0}
\end{equation}
that associates to $(X^-, \; X^0)$ the fixed point of \eqref{e:pert} is Frech\'{e}t differentiable and that its differential when $X^- (\tau) \equiv 0$ and $X^0 (\tau) \equiv \vec 0$
is the functional that associates to $(h^-, \, h^0) \in \mathcal Y$ the functions $U^-(\tau) \equiv \vec 0$, $U^0(\tau) \equiv \vec 0$. We then perform the linear change of variables which is the inverse of the change of variables introduced in Lemma \ref{l:lin}. In this way, we go back to the original variables. To simplify the notations, we still denote 
by $\Big(U^- (\tau), \; U^0 (\tau) \Big)$ the functions obtained applying the change of variables. 

To define the map that parameterizes the uniformly stable manifold we proceed as follows: the orbit $X^0(\tau)$ is fixed. For every $\underline X \in \mathbb R^{n_-}$, there exists a unique solution of \eqref{e:+:stac}. Also, in Section \ref{s:us:stable} we showed that the map
\begin{equation}
\label{e:p:ms}
    \underline X \to X^- (\tau) 
\end{equation}
is continuously differentiable in the sense of Frech\'{e}t. As a consequence, the map obtained composing \eqref{e:p:ms} and \eqref{e:p:fp:reg} is Frech\'{e}t differentiable.  Note that such a map associates to $\underline X^-$ the functions $(X^-,  \, U^-, \, U^0)$. The function $\phi$ that parameterizes the uniformly stable manifold is then defined by setting $$
     \phi (\underline X) = \Big( X^-(0) , \; X^0(0) + U^0(0) \Big).
$$
Thanks to the previous considerations, $\phi$ is continuously differentiable and the manifold is tangent to the stable space $ \big\{ (X^-, \, \vec 0): \; X^- \in \mathbb R^- \big\}$
at the origin. Also, estimate \eqref{e:p:exp:dec2} is a consequence of \eqref{e:us:es:u}. 

This concludes the proof of Theorem \ref{p:unif}.
\subsection{Uniformly stable manifolds}
\label{s:us:man}
Let $\mathcal V^0$ be a fixed center manifold for the equation 
\begin{equation}
\label{e:dec:e}
     \frac{d U}{d \tau} = F(U),
\end{equation}
which satisfies Hypotheses \ref{h:cutoff} and \ref{h:neg} introduced in Section \ref{s:hyp}. In Theorem \ref{p:unif} we consider a fixed orbit lying on $\mathcal V^0$ and we construct the uniformly stable manifold relative to that orbit. In this section we discuss what happens if, instead of having a single orbit, we have a whole invariant manifold. 

More precisely, let $\mathcal S_0$ be an invariant manifold for \eqref{e:us:theq:split} and assume that $\mathcal S_0$ is entirely contained in the center manifold $\{ X^- = \vec 0\}$. Also, denote by $\mathrm{S}^0$ the tangent space to $\mathcal S_0$ at the origin. Choosing a sufficiently small constant in Hypothesis \ref{h:cutoff}, we can assume that $\mathcal S_0$ is parameterized by $\mathrm{S}^0$. By construction, $\mathrm{S}^0$ is contained in $\{ X^- = \vec 0 \}$. Also, as in Section \ref{s:nont:zero} assume that $\mathcal Z_0= \{ (X^-, \, \vec 0, \, u_0): \zz = \vec 0\}$ is a manifold of zeroes for the function $f_0$ in \eqref{e:us:spli}.

As a consequence of Theorem \ref{p:unif}, we get the following result: 
\begin{pro}
\label{p:exp:d}
            Let Hypotheses \ref{h:cutoff} and \ref{h:neg} hold.  Let $\mathcal S_0$ be an invariant manifold for \eqref{e:us:theq:split} entirely contained in the center manifold $\{ X^- = \vec 0\}$. If the constant $\delta$ in Hypothesis \ref{h:cutoff} is sufficiently small, then the following holds.

           There exists a continuously differentiable manifold $\mathcal M^{us}_{\mathcal S_0}$ which is defined in the ball of radius $ \delta$ and center at the origin. Also, $\mathcal M^{us}_{\mathcal S_0}$ satisfies the following properties: 
           \begin{enumerate}
           \item $\mathcal M^{us}_{\mathcal S_0}$ is locally invariant for \eqref{e:us:theq:split}, meaning that if the initial datum lies on the manifold, then the solution $\big( (X^-(\tau), \, X^0(\tau) \big)$ of  \eqref{e:us:theq:split} lies on $\mathcal M^{us}_{\mathcal S_0}$  for $| \tau |$ sufficiently small. 
           \item $\mathcal M^{us}_{\mathcal S_0}$ is parameterized by $\mathrm{S}^0 \times V^-$ and it is tangent to this space at the origin. Here, $\mathrm{S}^0$ is the tangent space to $\mathcal S^0$ at the origin and $V^- = \{(X^-, \, \vec 0): \;  X^0 = \vec 0\}$. 
           \item Any orbit $Y(\tau)$ lying on $\mathcal M^{us}_{\mathcal S_0}$ can be decomposed as 
           \begin{equation}
           \label{e:dec:man:c}
                     Y(\tau) =  Y^0 (\tau) + Y^- (\tau) + Y^p (\tau), 
           \end{equation}
           where $Y^0 (\tau) = \big( \vec 0, \, \zz^0(\tau), \, u_0 (\tau) \big)$ is an orbit lying on $\mathcal S_0$. The component $Y^- (\tau) = \big( X^- (\tau), \, \vec 0, \, \vec 0 \big)$ lies on the stable manifold and the perturbation 
           term $Y^p (\tau)$ satisfies 
           \begin{equation}
           \label{e:p:exp:dec}
               |Y^p (\tau) | \leq C |  \zz^0 (0)| \, | Y^- (0)| e^{- c \tau / 4},
           \end{equation}
           for some positive constant $C$. In \eqref{e:p:exp:dec}, the constant $c >0$ is the same as in \eqref{e:us:Cmeno}.          
 \end{enumerate} 
\end{pro} 
 In the following we call $\mathcal M^{us}_{\mathcal S_0}$ the \emph{uniformly stable manifold relative to} $\mathcal S_0$.
\begin{proof}
Let $\big( \vec 0, \, X^0(\tau) \big)$ and $\big( X^-(\tau), \, \vec 0 \big)$ be two orbits of \eqref{e:us:theq:split} lying on the center manifold $\{ X^- = \vec 0\}$ and on the stable manifold respectively. We then have $X^0(\tau) \in \mathcal Y^0_{\delta}$, $X^-(\tau) \in \mathcal Y^-_{\delta}$, where the metric spaces $\mathcal Y^0_{\delta}$ and $\mathcal Y^-_{\delta}$ are defined by \eqref{e:st:y:d}. As in Section \ref{s:us:proof}, we use the notation $\mathcal Y = \mathcal Y^-_{\delta} \times \mathcal Y^0_{\delta}$. Consider the map
$$
    \Phi: \mathcal Y \to \mathcal Y \times \mathcal Y^p_{\delta 0}
$$   
which associates to $X^-(\tau)$ and $X^0(\tau) $ the function $\big( X^-(\tau), \, X^0(\tau), \, U^-(\tau), \, U^0 (\tau) \big)$, where $(U^-, \, U^0)$ is the perturbation term constructed in Section \ref{s:us:pert}. We recall that $\mathcal Y^p_{\delta 0}$ is the set obtained setting $a =0$ in \eqref{e:us:sp:ypd}. As shown in Section \ref{s:us:fre}, the map $\Phi$ is continuously differentiable in the sense of Frech\'{e}t. Also, let 
$$
     f^-:   \{ X^0 = \vec 0 \} \times  \{ X^- = \vec 0 \} \to \mathcal Y^-_{\delta}
$$
be the map that associates to $(\underline X^-, \, \underline \zz, \, \underline u^0 ) \in \mathcal Z_0$ the unique solution of the Cauchy problem  \eqref{e:+:stac}. We recall that in \eqref{e:+:stac} $\underline Y_0$ denotes $(\vec 0, \, \underline u_0)$. As shown in Section \ref{s:us:stable}, the map $f^-$ is continuously differentiable in the sense of Frech\'{e}t. Also, let 
$$
     f^0: \{ X^- = \vec 0 \}  \to \mathcal Y^0_{\delta}
$$
be the map that associates to $( \underline X^0, \, \vec 0 )$ the unique solution of the Cauchy problem \eqref{e:us:cau:x0}. The map $f^0$ is also continuously differentiable in the sense of Frech\'{e}t. Finally, fix a continuously differentiable map 
$$
     g^0: \mathrm{S^0} \to V^0 
$$   
parameterizing $\mathcal S_0$. Define the map
$$
     \psi:  \mathrm{S^0} \times  V^0 \to  \mathcal Y \times \mathcal Y^p_{\delta 0}
$$ 
setting
\begin{equation}
\label{e:f:psi}
    \psi( \underline X^0, \, \underline X^-) = \Phi \Bigg( f^- \Big(\underline X^-, \, g(\underline X^0) \Big), \, f^0 \circ g^0 (\underline X^0 ) \Bigg). 
\end{equation}
The map $\psi$ is then continuously differentiable in the sense of Frech\'{e}t. By construction, $  \psi( \underline X^0, \, \underline X^-)$ is an element in the form 
$\big( X^-(\tau), \, X^0(\tau), \, U^-(\tau), \, U^0(\tau) \big)$ and, setting 
$$
     Y(\tau) = \Big( X^-(\tau) + U^-(\tau), \; X^0(\tau)+ U^0(\tau) \Big),
$$
we get that $Y(\tau)$ can be decomposed as in \eqref{e:dec:man:c}. Also,  the perturbation term $\big( U^0, \, U^- \big)$ automatically satisfies \eqref{e:p:exp:dec}. We then define the map
$$
    \psi_0:   
    \mathrm{S^0} \times  V^- \to \mathbb R^{n_c + n_-}
$$ 
parameterizing $\mathcal M^{us}_{\mathcal S_0}$ by setting 
$$
     \psi_0 ( \underline X^0, \, \underline X^-) = \Big( X^-(0) , \; X^0(0) + U^0( 0) \Big) = Y(0),
$$
where $X^0(\tau)$, $X^-(\tau)$ and $U^0(\tau)$ are given by \eqref{e:f:psi}. 

The map $\psi_0$ is continuously differentiable, being the composition of maps that are continuously differentiable in the sense of Frech\'{e}t.  Also, by construction the manifold $\mathcal M^{us}_{\mathcal S_0}$ is invariant for \eqref{e:us:theq:split}. To prove that the manifold $\mathcal M^{us}_{\mathcal S_0}$ is tangent to $\mathrm S^0 \times V^-$ at the origin it is enough to observe that the Frech\'{e}t differential of $f^-$ at $\underline X^- = \vec 0$ is the functional $\underline h_- \mapsto e^{\bar A^- \tau} \underline h^-$, while the Frech\'{e}t differential of $f^0$ at $\underline X^0 = \vec 0$ is the functional $\underline h_0 \mapsto e^{\bar A^0 \tau} \underline h^0$.   
\end{proof}

\section{Invariant manifolds for a singular ODE}
\label{s:invman}
In Section \ref{s:invman} we extend to the general case the considerations introduced in Section \ref{sus:intro_toy} in the case of a toy model. In doing this, we apply the results obtained in Section \ref{s:us} to study the singular ordinary differential equation
\begin{equation}
\label{e:man:sing}
          \frac{d U}{dt} = \frac{1}{\zz (U)} F(U) .
\end{equation}
Actually, most of the time we focus on system 
\begin{equation}
\label{e:m:orig}
          \frac{d U}{ d \tau} = F(U) .
\end{equation}
We discuss several situations where \eqref{e:man:sing} and \eqref{e:m:orig} are equivalent, namely the Cauchy problem 
$$
\left\{
\begin{array}{ll}
           d \tau  / d t = 1 / \zz [ U(t) ] \\
           \tau (0) = 0
\end{array}
\right.
$$
defines a continuously differentiable diffeomorphism $\tau: [0, \, + \infty [ \to [0, \, + \infty [ $.

In Section \ref{s:invman} we exploit all Hypotheses \ref{h:pos} $\dots$ \ref{h:slow}. Also, we rely on Proposition \ref{p:change}, whose proof is given in Section \ref{sus:m:change}. Before stating it, we have to introduce some notations. Let     
 $N$ denote the dimension of $U$. Also, $n_-$ is the number of eigenvalues of 
          $D F (\vec 0)$ with strictly negative real part, while $(n_0 + 1)$ is the number of eigenvalues of 
          $D F (\vec 0)$ with zero real part. Each eigenvalue is counted according to its multiplicity. Thanks to Hypothesis \ref{h:neg}, $N= n_-+ n_0+1$.       
\begin{pro}
\label{p:change}
          Let Hypotheses \ref{h:pos} $\dots$ \ref{h:slow}  hold. If the constant $\delta$ in Hypothesis \ref{h:cutoff} is sufficiently small, then in the ball with radius $\delta$ and center at the origin we can define a continuously differentiable diffeomorphism $\uu$ satisfying the following properties. 
          Write $\uu (U) = \bar U$ as a column vector: 
          $$ 
               \bar U = 
               \left(
               \begin{array}{ccc}
                           \zeta \\
                           u_0 \\
                           u_- \\
               \end{array}
              \right),
          $$
          where $ \zeta \in \mathbb R$, $u_0 \in \mathbb R^{n_0}$ and $ u_- \in \mathbb R^{n_-}$. If $U$ satisfies \eqref{e:m:orig} then $\bar U $ satisfies
	 \begin{equation}
	   \label{e:m:equation}
	   \left\{
     \begin{array}{lllll}
                 \displaystyle{ 
		 d  \zeta/ d \tau=        G_{10} ( \zeta,  \,  u_0 )  u_0 \zeta^2  +  G_{1 - } ( \zeta,  \,  u_0, \,  u_-)  u_-  \zeta } \\
		 \\
      \displaystyle{ d  u_0 / d \tau = \Big\{     G_{0 1} ( \zeta, \,  u_0)  + \big[   G_{0-} ( \zeta, \,  u_0, \,  u_-) -      G_{0-} ( \zeta, \,  u_0, \, \vec 0) \big]  \Big\}  \zeta  u_0 } \\
      \\
           \displaystyle{d  u_-/ d \tau =   G_s  ( \zeta, \,  u_0, \,  u_-)  u_-  } \\
     \end{array}
     \right.
     \end{equation}     
     In the previous expression, $ G_{1 0}$ is a row vector belonging to $\mathbb R^{n_0}$, $G_{1 -}$ is a row vector in $\mathbb R^{n_-}$, 
     the matrices $G_{01}$ and $ G_{0-}$ belong to $\mathbb M^{n_0 \times n_0}$ and the matrix $ G_s$ belongs to $\mathbb M^{n_- \times n_-}$.  
    
     A center manifold of system \eqref{e:m:equation} is the subspace $\{ (\zz, \, u_0, \, \vec 0): \; u_- = \vec 0 \}$, the stable manifold is the subspace 
     $\{ (0 , \, \vec 0, \, u_-): \; \zz =0, \; u_0= \vec 0 \}$. Let $\mathcal M^{us}_E$ be the uniformly stable manifold relative to the manifold ${E= \{  (\zz, \, \vec 0, \, \vec 0): \; u_0 = \vec 0, \; u_- = \vec 0 \}}$, which is entirely constituted by equilibria. Then ${\mathcal M^{us}_E = \{ (\zz, \, \vec 0, \, u_-): \; u_0 = \vec 0 \}.}$
      \end{pro}
In the statement of Proposition \ref{p:change} by \emph{uniformly stable manifold relative to} $E$ we mean the manifold defined by Proposition \ref{p:exp:d}. Also, note that 
 by construction all the eigenvalues of the matrix $ G_s(0, \, \vec 0, \, \vec 0)$ have strictly negative real part.
\subsection{Slow and fast dynamics}
\label{sus:slowfast}
Let $E$ denote, as before, the manifold of equilibria $\{  (\zz, \, \vec 0, \, \vec 0): \; u_0 = \vec 0, \; u_- = \vec 0 \}$.
\begin{say}
\label{say:slowfast}
         A manifold of slow dynamics is a center manifold of \eqref{e:m:equation}. 
      In the following we fix the manifold of the slow dynamics $\{ u_- = \vec 0 \}$ and we denote it by $\mathcal M^0$.  
      
The manifold of fast dynamics of system \eqref{e:m:equation} is the uniformly stable manifold relative to $E$, namely the subspace $\{ u_0 = \vec 0  \}$.  
\end{say}
Note that both these manifolds are invariant for system \eqref{e:m:equation}. Also, for every point $(\underline \zz, \, \vec 0, \, \underline u_-)$ belonging to the manifold of fast dynamics, denote by $\big(\zz (\tau), \, \vec 0, \, u_-(\tau) \big)$ the solution of \eqref{e:m:equation} such that   
$$ 
  \big(   \zz (0) , \, \vec 0, \, u_-(0) \big)=(\underline \zz, \, \vec 0, \, \underline u_-).
$$
Combining \eqref{e:dec:man:c} and \eqref{e:p:exp:dec} we get that this solution decays exponentially fast to an equilibrium point.  Namely, there exists $ \zz_{\infty}$ such that 
$$
   \lim_{\tau \to + \infty} e^{c \tau / 4}| u_- (\tau) |  = 0 =   \lim_{\tau \to + \infty} e^{c \tau / 4}| \zz (\tau) - \zz_{\infty} |,
$$
where the positive constant $c$ satisfies $Re \lambda < -c$ for every $\lambda$ eigenvalue of $G_s (0, \, \vec 0, \, \vec 0)$.

Consider system \eqref{e:m:equation} reduced on the manifold of slow dynamics: 
\begin{equation}
\label{e:m:reduced}
  \left\{
     \begin{array}{lll}
                 \displaystyle{ 
		 d \zeta/ d \tau= 
		  \zeta^2 G_{10} ( \zeta,  \,  u_0 )  u_0}  \\
                  \displaystyle{ d  u_0 / d \tau =   G_{0 1} ( \zeta, \,  u_0) } u_0 \zeta  \\
                u_- \equiv 0 \\
     \end{array}
     \right.
\end{equation}   
If one goes back to the original variable $t$ obtains   
\begin{equation}
\label{e:m:t}
  \left\{
     \begin{array}{lllll}
                 \displaystyle{ 
		 d \zeta/ d t = 
		 \zeta  G_{10} ( \zeta,  \,  u_0 )  u_0 } \\
                  \displaystyle{ d  u_0 / d t =   G_{0 1} ( \zeta, \,  u_0)  } u_0 \\
                u_- \equiv 0, \\
     \end{array}
     \right.
\end{equation} 
namely an equation with no singularity. Note that \eqref{e:m:reduced} and \eqref{e:m:t} are equivalent. Indeed, by the uniqueness of the solution of a Cauchy problem associated to \eqref{e:m:t}, the following holds. If $\zz (0) > 0$ then $\zz (t) > 0$ for every $t$. Thus, the Cauchy problem 
\begin{equation}
\label{e:m:cauchy}
  \left\{
  \begin{array}{lll}
      \displaystyle{\frac{d \tau}{d t } = \frac{1}{\zz (t)} }\\
          \\
        \tau (0) =0 \\
  \end{array}
  \right.
\end{equation}
admits a global solution $\tau: [0, \, + \infty [ \to [0 , \, + \infty [ $ whose derivative is always different from $0$. Thus, $\tau (t)$ defines a change of variables and \eqref{e:m:reduced} is equivalent to \eqref{e:m:t}. 

One of our original goals is to study the solutions of 
$$
  \frac{d U}{d t} = \frac{\phi_s (U)}{\zeta (U)} + \phi_{ns}(U)
$$
lying on a center manifold. Let $\mathcal M^{00}$ be a center manifold for \eqref{e:m:t} around the equilibrium point $(0, \, \vec 0, \, \vec 0)$. Then $\mathcal M^{00}$ is a center manifold for 
  \begin{equation}
  \label{e:m:cnozeta}
  \left\{
     \begin{array}{lllll}
                 \displaystyle{\frac{ 
		 d  \zeta}{d t}=        \zeta G_{10} ( \zeta,  \,  u_0 )  u_0+    G_{1 - } ( \zeta,  \,  u_0, \,  u_-)  u_-} \\
		 \\
      \displaystyle{ \frac{d  u_0 }{d t} = \Big\{     G_{0 1} ( \zeta, \,  u_0)  + \big[   G_{01} ( \zeta, \,  u_0, \,  u_-) -      G_{0-} ( \zeta, \,  u_0, \, \vec 0) \big]  \Big\}   u_0 } \\
      \\
           \displaystyle{\frac{d  u_-}{d t } =   \frac{1}{\zeta }G_s  ( \zeta, \,  u_0, \,  u_-)  u_-  } \\
     \end{array}
     \right.
 \end{equation}
 We collect these results in the following
 \begin{teo}
 \label{t:center}
        Assume that Hypotheses \ref{h:pos} $\dots$ \ref{h:slow} are satisfied. 
        There exists an invariant  center manifold $\mathcal M^{00}$ for system 
       \eqref{e:m:cnozeta} around the equilibrium point $(0, \, \vec 0, \, \vec 0)$ which is contained in the manifold of the slow dynamics. In particular, equation \eqref{e:m:cnozeta} restricted to $\mathcal M^{00}$ is non singular and every solution satisfies the following property: if $\zz (0) > 0$, then $\zz(t) > 0$ for every $t$.  
 \end{teo}
 \begin{rem}
 \label{r:slow}
           Hypothesis \ref{h:slow} ensures that the manifold $\{ U: \; \zz (U) =0 \}$ is invariant with respect to the slow dynamics.  This hypothesis is not necessary to define an invariant center manifold $\mathcal M^{00}$ contained in the manifold of the slow dynamics. However, it is necessary if we want that \eqref{e:m:reduced} is equivalent to \eqref{e:m:t}, namely that the change of variables defined by \eqref{e:m:cauchy} is well defined.  
         To see this, we can proceed as follows. 
           
           Consider the equation 
           $$
               \frac{d U}{ d \tau} = F(U).
           $$
           Assume that one proceeds as in the proof of Lemma \ref{l:ch} and exploits Hypotheses \ref{h:pos} $\dots$ \ref{h:fast} but  does \emph{not} exploit Hypothesis \ref{h:slow}. The system one eventually gets, restricted on the subspace $\{ u_- = \vec 0 \}$, is  
\begin{equation}
\label{e:r:slow}
  \left\{
     \begin{array}{lll}
                 \displaystyle{ 
		 d \zeta/ d \tau= 
		  \zeta g_{1} ( \zeta,  \,  u_0, \, \vec 0)}  \\
                  \displaystyle{ d  u_0 / d \tau =   G_{0 1} ( \zeta, \,  u_0) } u_0 \zeta  \\
                                  u_- \equiv 0 \\
     \end{array}
     \right.
\end{equation}              
where $g_1$ is the same function as in \eqref{e:m:system1} and satisfies 
$$
    g_1 (z, \, \vec 0, \, \vec 0) = 0 \quad \forall \, z .
$$           
 
Going back to the original variable $t$, \eqref{e:r:slow} becomes 
 \begin{equation}
\label{e:m:r}
  \left\{
     \begin{array}{lll}
                 \displaystyle{ 
		 d \zeta/ d \tau= 
		  g_{1} ( \zeta,  \,  u_0, \, \vec 0)}  \\
                  \displaystyle{ d  u_0 / d \tau =   G_{0 1} ( \zeta, \,  u_0) } u_0   \\
                u_- \equiv 0 \\
     \end{array}
     \right.
\end{equation}
Thus, even if we do \emph{not} assume Hypothesis \ref{h:slow}, the equation
$$
   \frac{d U }{dt} = F(U)
$$
restricted on the manifold of the slow dynamics $\{ u_- = \vec 0 \}$ is non singular. Also, one can define an invariant center manifold $\mathcal M^{00}$ which contains only slow dynamics.

Note, however, that if Hypothesis \ref{h:slow} is not satisfied it may happen that for a solution $U$ lying on $\mathcal M^{00}$ $ \zz \big(U (0) \big) >0$ but $\zz (U)$ touches $0$ in finite time. An example is the following. 

Consider the equation 
$$
\left\{
\begin{array}{lll}
            d u_1 / dt = - u_2 \\
            d u_2 / dt = u_2^2 ( 1 - u_2) \\
            d u_3 / dt = - u_3 / u_1
\end{array}
\right.
$$
and set 
$$
   \zz (U) = u_1 \qquad F(U) = \Bigg( -u_1 u_2, \; u_1 u_2^2 ( 1 - u_2), \, -u_3  \Bigg)^T.
$$
Then Hypotheses \ref{h:pos}, \ref{h:neg} $\dots$ \ref{h:fast} are satisfied, but Hypothesis \ref{h:slow} is violated. The manifold of slow dynamics is the subspace $\{ u_3 =0 \}$ and it coincides with the center manifold $\mathcal M^{00}$. 
Restrict to this subspace and consider the equation
$$
     d u_2 / dt = u_2^2 ( 1 - u_2).
$$
If $0 < u_2 (0) <1 $, then $0 < u_2 (t) <  1$ for every $t$. Also, $d u_2 / dt > 0$ for every $t$ and hence $u_2 (0) < u_2 (t) < 1 $ for every $t$. Since 
$$
       d u_1 / dt = - u_2 ,
$$
then by a comparison argument $u_1 (t ) \leq u_1 (0) - u_2 (0) t $ for every $t >0$. In other words, if $u_1 (0) >0$ then $u_1(t)$ attains the value $0$ for some $t \leq u_1 (0) / u_2(0)$. 
\end{rem}
 
\subsection{Applications of the uniformly stable manifold to the analysis of a singular ordinary differential equation}
\label{s:appl}
We first recall a preliminary result we need in the following
\begin{lem}
\label{l:blow:up}
          Let $ \zz(\tau)$ be a real valued, continuous and bounded function satisfying $\zz (\tau) >0$ for every $\tau \in [0, \, + \infty[$. Let $t(\tau)$ be the maximal solution of the forward Cauchy problem 
          \begin{equation}
          \label{e:app:c:f}
           \left\{
           \begin{array}{ll}
                   d t / d \tau = \zz (\tau) \\
                   t (0) = 0 \\
           \end{array}
           \right.
          \end{equation}
          Then $t(\tau)$ is defined on the whole interval $[0, \, + \infty[$. Also, the following statements are equivalent:
          \begin{enumerate}
          \item $t(\tau)$ is a continuously differentiable diffeomorphism $t:[0, \, + \infty[ \to [0, \, + \infty[$. 
          \item $\displaystyle{\int_0^{+\infty} \zz (\tau) d \tau}= + \infty$.
          \end{enumerate}
\end{lem}
Condition $2$ guarantees, in particular, that the inverse map $\tau(t)$ is defined on the whole interval $[0, \, + \infty[$ and that it is continuously differentiable there. Also, note that $\zz (t) = \zz \big( \tau (t) \big) $ is automatically strictly bigger than $0$ for every $t$. 

Before stating the most important result in this section we need to introduce some notations.
As before, $c>0$ denotes a positive constant satisfying $Re \lambda < -c $ for any $\lambda$ which is either an eigenvalue of $G_s (0, \, \vec 0, \, \vec0)$ or an eigenvalue with strictly negative real part of  
of $G_{01} (0, \, \vec 0)$. We denote by 
$V^{0-}$ the subspace
$$
     V^{0-} =\{(0, \, \vec \xi, \, \vec 0)  \},
$$
where $\vec \xi \in \mathbb R^{n_0}$ belongs to the eigenspace of $G_{10}(0, \, \vec 0)$ associated to the eigenvalues with strictly negative real part. Also, 
$$
     V^{00-} =\{(0, \, \vec \xi, \, \vec 0)  \},
$$
where $\vec \xi \in \mathbb R^{n_0}$ belongs to the eigenspace of $G_{10}(0, \, \vec 0)$ associated to the eigenvalues with non positive real part. Clearly, $V^{0-} \subseteq V^{00-}$. With $V^{--}$ we denote the stable manifold:
$$
     V^{--} =\{(0, \, \vec 0, \,  u_-): \; u_- \in \mathbb R^{n_-}  \},
$$
Finally, as in Section \ref{sus:slowfast} we denote by $E$ the manifold of equilibria  $\{ (\zeta, \, \vec 0, \, \vec 0): \; \zz \in \mathbb R \}$.

The most important result in this section is the following:
\begin{teo}
\label{t:app:man}
           Let Hypotheses \ref{h:pos} $\dots$ \ref{h:slow} hold. If the constant $\delta$ in Hypothesis \ref{h:cutoff} is sufficiently small, then in the ball with radius $\delta$ and center at the origin one can define two manifolds, $\mathcal M^s$ and $\mathcal M^{cs}$, satisfying the following properties:
           \begin{enumerate}
           \item both $\mathcal M^s$ and $\mathcal M^{cs}$ are locally invariant for \eqref{e:m:equation}, namely: if the initial datum lies on the manifold, then the solution $\big( (\zz (\tau), \, u_0(\tau), \, u_-(\tau)  \big)$ of  \eqref{e:m:equation} also lies on the manifold  for $| \tau |$ sufficiently small. 
           \item $\mathcal M^s$ is contained in $\mathcal M^{cs}$. 
           \item $\mathcal M^s$ is parameterized by $E \oplus V^{0-} \oplus V^{--}$ and it is tangent to this subspace at the origin. Also,  $\mathcal M^{cs}$ is parameterized by $E \oplus V^{00-}  \oplus V^{--}$ and it is tangent to this subspace at the origin. 
           \item let $U(\tau) = \Big( \zz (\tau), \; u_0 (\tau), \, u_- (\tau) \Big)$ be an orbit lying either on $\mathcal M^s$ or on $\mathcal M^{cs}$ and satisfying $\zz (0) >0$. 
           Then the maximal solution of the forward Cauchy problem 
            \begin{equation}
          \label{e:app:c:f2}
           \left\{
           \begin{array}{ll}
                   d t / d \tau = \zz (\tau) \\
                   t (0) = 0 \\
           \end{array}
           \right.
          \end{equation}
           defines a continuously differentiable diffeomorphism $t:[0, \, + \infty[ \to [0, \, + \infty[$. Let $\tau(t)$ denotes its inverse. Then the function $U (t) = U \big( \tau (t) \big)$ is a solution of \eqref{e:man:sing} and satisfies $\zz (t) >0 $ for every $ t \ge 0$.   
           \item any orbit lying on $\mathcal M^s$ can be decomposed as 
           \begin{equation}
           \label{e:app:t:s}
                U(\tau) = U^- (\tau) + U^{sl} (\tau) + U^p (\tau),
           \end{equation}
           where $U^-(\tau) $ satisfies 
           \begin{equation}
           \label{e:app:t:m}
                     |U^- (\tau)| \leq k_- e^{- c \tau /2} |U^-(0)|
           \end{equation} 
           for a suitable constant $k_-$.  Conversely, the component $U^{sl} (\tau) = \big( \zz^{sl} (\tau), \, u_0^{sl} (\tau), \, \vec 0 \big)$ lies on the manifold of the slow dynamics. Also, if we use the variable $t$ defined as the maximal solution of the Cauchy problem \eqref{e:app:c:f2}, 
           we have that the following property is satisfied. Denote by $\zz$ and $u_0$ the first and the second component of $U$ respectively. Then there exists a point $(\zz_{\infty}, \, \vec 0)$ such that 
           \begin{equation}
           \label{e:app:t:d}
                \lim_{t \to + \infty} \big( |\zz (t) - \zz_{\infty} | + | u_0 (t) | \big) e^{ c t / 2}  =0.
           \end{equation}

           Finally, the perturbation term is small in the sense that 
           \begin{equation}
           \label{e:us:small:p}
                |U^p (\tau) | \leq k_p | \zz^{sl} (0) | |U^-(0) | e^{- c \tau / 4}
           \end{equation}
           for a suitable constant $k_p >0$. 
           \item any orbit $U (\tau)$ lying on $\mathcal M^{cs}$ can be decomposed as 
           \begin{equation}
           \label{e:app:t:cs}
                U(\tau) = U^- (\tau) + U^{sl} (\tau) + U^p (\tau),
           \end{equation}
           where $U^-$ and $U^p$ satisfy $|U^- (\tau)| \leq k_- e^{- c \tau /2} |U^-(0)|$ and $|U^p (\tau) | \leq k_p | \zz^{sl} (0) | |U^-(0) | e^{- c \tau / 4}$ respectively. Here $k_-$ and   $k_p$ denote the same constants as in \eqref{e:app:t:m} and \eqref{e:us:small:p}. The component $U^{sl}(\tau)= \big( \zz^{sl} (\tau), \, u^{sl} (\tau), \, \vec 0 \big)$ lies on the manifold of the slow dynamics. More precisely, the following holds. Consider the maximal solution of the Cauchy problem
           $$
              \left\{
              \begin{array}{ll}
                   d t / d \tau = \zz^{sl} (\tau) \\
                   t (0) = 0 \\
     \end{array}
     \right.
     $$
    and set
          $\zz^{sl}(t) = \zz^{sl} \big( \tau (t) \big)$ and $u^{sl}(t) = u^{sl} \big( \tau (t) \big)$.  Then $\big( \zz^{sl}(t), \, u^{sl}(t)  \big)$ is a solution lying on a center-stable manifold of 
          $$
              \left\{
     \begin{array}{lllll}
                 \displaystyle{ 
		 d \zeta/ d t = 
		 \zeta  G_{10} ( \zeta,  \,  u_0 )  u_0 } \\
                  \displaystyle{ d  u_0 / d t =   G_{0 1} ( \zeta, \,  u_0)  } u_0 \\
                u_- \equiv 0, \\
     \end{array}
     \right.
          $$
          \end{enumerate}
           \end{teo}
           Note that, strictly speaking, in \eqref{e:app:t:s} and in  \eqref{e:app:t:cs} the component $U^-$ does \emph{not} lie on the manifold of the fast dynamics. Indeed, as we will see in the proof, $U^-$ is a solution of \eqref{e:+:stac} and hence does not lie on $\{ (0, \, \vec 0, \, u_- ) \}$. However, loosely speaking it can be regarded as a fast dynamic because of its exponential decay.  
\begin{proof}
We first define $\mathcal M^{s}$.

Consider system \eqref{e:m:equation} restricted on the manifold of the slow dynamics. Thanks to the analysis in Section \ref{sus:slowfast} the variables $t$ and $\tau$ are then equivalent. Using the variable $t$, we get  
\begin{equation}
\label{e:f:t}
  \left\{
     \begin{array}{lll}
                 \displaystyle{ 
		 d \zeta/ d t = 
		 \zeta  G_{10} ( \zeta,  \,  u_0 )  u_0 } \\
                  \displaystyle{ d  u_0 / d t =   G_{0 1} ( \zeta, \,  u_0)  } u_0 \\
                u_- \equiv 0, \\
     \end{array}
     \right.
\end{equation} 
The manifold $E= \{ (\zeta, \, \vec 0, \, \vec 0): \; \zz \in \mathbb R \}$ is then entirely constituted by equilibria. Applying Proposition \ref{p:exp:d} to system \eqref{e:f:t} with $\mathcal S_0 = E$, we then obtain $M^{us}_E$, the uniformly stable manifold relative to $E$, which is parameterized by $E \oplus V^{0-}$. Note that so far we have used only the variable $t$: $M^{us}_E$ is a uniformly stable manifold for \eqref{e:f:t} with respect to the variable $t$ and by construction it is included in $\{ u_- = \vec 0 \}$, a center manifold for \eqref{e:m:equation} with respect to the variable $\tau$. The manifold $\mathcal M^{s}$ is then obtained exploiting the variable $\tau$ and applying Proposition \ref{p:exp:d} to system \eqref{e:m:equation} with $\mathcal S_0 = M^{us}_E$. Also, the set
$$ 
      \mathcal Z_0 = \{ (0, \, u_0, \, u_-): \; u_0 \in \mathbb R^{n_0}, \; u_- \in \mathbb R^{n_-} \}. 
$$   
satisfies \eqref{e:us:zero}. 
Properties 1, 3 and estimates \eqref{e:app:t:m} and \eqref{e:us:small:p} in the statement of Theorem \ref{t:app:man} are then automatically satisfied, so we are left to prove estimate \eqref{e:app:t:d} and property 4. 

To show that estimate \eqref{e:app:t:d} holds we apply Lemma \ref{l:blow:up}. Thanks to \eqref{e:app:t:s}, 
$$
     \zz (\tau) = \zz^{sl}(\tau) + \zz^p (\tau),
$$
where $U^{sl} (\tau) = \big( \zz^{sl}(\tau), \, u_0^{sl} (\tau), \, \vec 0 \big)$ lies on the manifold of the slow dynamics and $\zz^p$ is the first component of the perturbation term $U^p$. 
Let $\tilde t$ be defined as the maximal solution of 
           \begin{equation*}         
           \left\{
           \begin{array}{ll}
                   d \tilde t  / d \tau = \zz^{sl} (\tau) \\
                   \tilde t (0) = 0, \\
           \end{array}
           \right.
           \end{equation*}
Then there exits $(\zz_{\infty}, \, \vec 0)$ such that 
        \begin{equation}
        \label{e:app:t:d2}
                \lim_{\tilde t \to + \infty} \big( |\zz^{sl} (\tilde t) - \zz_{\infty} | + | u^{sl}_0 (\tilde t) | \big) e^{ c \tilde t / 2}  =0. 
      \end{equation}
Since $|\zz^p (\tau)| \leq k_p \delta^2 e^{- c \tau /4}$, then for every $\tau$
$$
   |   \tilde t (\tau) - t (\tau) | \leq \unpo \delta^2
$$
where $t (\tau)$ is defined by \eqref{e:app:c:f2}. Since also $| u^{sl}_0 (\tau) - u^{\sl}_0 (\tau) | \leq k_p \delta^2 e^{- c \tau /4}$, we conclude that \eqref{e:app:t:d} implies \eqref{e:app:t:d2}. 

Concering the proof of property 4, we apply Lemma \ref{l:blow:up}. Since $U^{sl} (\tau) = \big( \zz^{sl}(\tau), \, u_0^{sl} (\tau), \, \vec 0 \big)$ lies on the manifold of the slow dynamics, then by the analysis in Section \ref{sus:slowfast} it satisfies condition 1 in the statement of Lemma \ref{l:blow:up} and hence 
$$
     \int_0^{+ \infty} \zz^{sl} (\tau) d \tau = + \infty. 
$$
Since $|\zz^p (\tau)| \leq \delta^2 e^{- c \tau /4}$, then 
$$
     \int_0^{+ \infty} \big( \zz^{sl} + \zz^p \Big)  (\tau) d \tau = + \infty. 
$$
Applying again Lemma \ref{l:blow:up} we get property 4. 

To define the manifold $\mathcal M^{cs}$ we proceed as follows.  Consider $M^{cs}$, a center-stable manifold for \eqref{e:f:t}. This manifold is parameterized by $E  \oplus V^{00-}$ and it is tangent to this space at the origin. The manifold $\mathcal M^{cs}$ is defined applying Proposition \ref{p:exp:d} to system \eqref{e:m:equation} with $\mathcal S_0 = M^{cs}$ and exploiting the presence of the set $\mathcal Z_0 = \{ (0, \, u_0, \, u_-): \; u_0 \in \mathbb R^{n_0}, \; u_- \in \mathbb R^{n_-} \}$ satisfying \eqref{e:us:zero}. Proceeding as before one gets that properties 1, 3, 4 and 6 are satisfied. 

To verify property 2, we first observe that $M^{us}_E \subseteq M^{cs}$. To obtain $\mathcal M^{s}$ and $\mathcal M^{cs}$ we applied Proposition~\ref{p:exp:d}  to $\mathcal S_0= M^{us}_E$ and $\mathcal S_0=M^{cs}$ respectively.  Going back to the proof of Proposition \ref{p:exp:d}  one can notice that the way we constructed the uniformly stable manifold with respect to $\mathcal S_0$ is we associated to any orbit lying on $\mathcal S_0$ the manifold constructed in Theorem \ref{p:unif}. Thus the inclusion $M^{us}_E \subseteq M^{cs}$  
has as a consequence the inclusion  $\mathcal M^{s} \subseteq \mathcal M^{cs}$.
\end{proof}

\subsection{Proof of Proposition \ref{p:change}}
\label{sus:m:change}

\subsubsection{A preliminary result}
\label{sus:pre}
Before proving Proposition \ref{p:change}, we have to introduce a preliminary result, Lemma \ref{l:ch}. 

Let $\Upsilon$ be a continuously differentiable local diffeomorphism.  To simplify the exposition, we also assume that $ \Upsilon (\vec 0) = \vec 0$.  
Let $\bar U: = \Upsilon (U)$ and 
\begin{equation}
\label{e:m:barF}
    \bar F (\bar U) : = D  \Upsilon \big( \bar \Upsilon^{-1 } (\bar U) \big) F \big( \Upsilon^{-1} (\bar U) \big)
\end{equation}
If the function $U (\tau)$ satisfies \eqref{e:m:orig}, then $\bar U(\tau)$ solves 
\begin{equation}
\label{e:m:bar}
          \frac{d \bar U}{d \tau } = \bar F (\bar U). 
\end{equation}
Also, given a real valued function $\zeta (\bar U)$, let
\begin{equation}
\label{e:m:barzeta}
     \bar \zeta (\bar U) : = \zeta \big[ \Upsilon^{-1} (\bar U) \big].
\end{equation}
By direct check, one can verify that the following holds true. 
\begin{lem}
\label{l:ch}
           Assume that Hypotheses \ref{h:pos}, \ref{h:neg} $\dots$ \ref{h:slow} are satisfied by $F$ and $\zeta$. Also, assume that Hypothesis \ref{h:cutoff} is satisfied for some $\delta$. Then Hypotheses \ref{h:pos}, \ref{h:neg} $\dots$ \ref{h:slow}  are verified by $\bar F$ and $\bar \zeta$ and there exists $\bar \delta$, possibly smaller than $\delta$,  such that Hypothesis \ref{h:cutoff} is as well satisfied.  
\end{lem}
\subsubsection{Proof of Proposition \ref{p:change}: first part}
\label{sus:first part}
We are now ready to prove Proposition \ref{p:change}. The proof actually relies on standard techniques, but we give it for completeness.  We proceed in several steps. 
\begin{itemize}
\item {\it Step 1:}
let $U= (u_1 \dots u_N)^T$ be the components of $U$. Thanks to Hypothesis \ref{h:sur}, $\nabla \zz (\vec 0) \neq \vec 0$. Just to fix the ideas, we can assume 
$$
    \frac{\partial \zeta}{\partial u_1 } (\vec 0) \neq 0. 
$$
By a smooth local change of variables we can assume that $\zeta(U) = u_1$. Thanks to Lemma \ref{l:ch}, Hypotheses \ref{h:pos} $\dots$ \ref{h:slow} are satisfied by the ODE written using the  new variable. To simplify the exposition, we write $U$ and $\zz$ instead of $\bar U$ and $\bar \zz$. 
\item {\it Step 2:} thanks to Hypothesis \ref{h:tras}, there exists a manifold $\mathcal M^{eq}$ which is 
entirely constituted by equilibria and which is transversal to the manifold $\mathcal S$, namely to $\{ u_1 = 0 \}$.  \emph{Via} a smooth local 
change of variables we can assume that the one-dimensional subspace 
\begin{equation}
\label{e:m:e}
         E : = \{ \bar u_2 = \dots = \bar u_N = 0 \}
\end{equation}
 is entirely contained in $\mathcal{M}^{eq}$. Hypotheses \ref{h:pos} $\dots$ \ref{h:slow} are satisfied in the new variables thanks to Lemma \ref{l:ch}.  
\item {\it Step 3:}  let $E$ be as in \eqref{e:m:e} and denote by $V^c$ the eigenspace of $DF(\vec 0)$ associated to eigenvalues with $0$ real part.  Also, let $V^{--}$ be the eigenspace associated to eigenvalues with strictly negative real part. The dimension of $V^c$ and of $V^{--}$ is $n_0 +1$ and $n_-$ respectively. Thanks to Hypothesis \ref{h:neg}, $N = n_0 +1 + n_-$. 
The vector $(1, \, 0 \dots 0)$ belongs to $V^c$ because $E \subseteq V^c$. Also, we can assume, \emph{via} a linear change of variables, that 
$$
     V^c = \{ u_{n_0 +2  } = \dots u_N =0  \} \qquad  V^s = \{ \zz =0, \, u_{2} = \dots u_{n_0 + 1} = 0  \}.  
$$
Fix a center manifold $\mathcal M^c$ for system 
\begin{equation}
\label{e:u:s}
    \frac{d U }{d \tau} = F(U)
\end{equation}
around the equilibrium point $\vec 0$: $\mathcal M^c$ is parameterized by $V^c$ and it is tangent to this space at the origin $\vec 0$. Also, let $\mathcal M^{us}_E$ be the uniformly stable manifold of system \eqref{e:u:s} relative to the manifold of equilibria $E$ defined by \eqref{e:m:e}: this manifold is paramerized by $V^s \oplus E$ and it is tangent to this space at the origin.  By a local smooth change of variables we can assume that actually
$$
     \mathcal M^c = \{ u_{n_0 +2  } = \dots u_N =0  \}  \qquad \mathcal M^{us}_E =   \{  u_{2} = \dots u_{n_0 + 1} = 0  \}. 
$$
Note that the Hypotheses \ref{h:pos} $\dots$ \ref{h:slow} are satisfied because of Lemma \ref{l:ch}.  
\item {\it Step 4:} consider the following decomposition:
\begin{equation}
\label{e:m:dec}
         U =
         \left(
         \begin{array}{ccc}
                     \zeta \\
                     u_0 \\
                     u_-
         \end{array}
         \right)
         \qquad 
         F(U) =
         \left(
         \begin{array}{ccc}
                     f_1 (\zeta, \, u_0, \, u_-) \\
                     F_0   (\zeta, \, u_0, \, u_-)  \\
                     F_-  (\zeta, \, u_0, \, u_-), \\
         \end{array}
         \right)
\end{equation}
where $\zeta, \; f_1 \in \mathbb R$, $u_0, \; F_0 \in \mathbb R^{n_0}$ and $u_-, \; F_- \in \mathbb R^{n_-}$. In the new coordinates, the center manifold $\mathcal M^c$ is the subspace $\{ u_- = \vec 0 \}$ and the uniformly stable manifold $\mathcal M^{us}_E$ is $\{ u_0 = \vec 0 \}$.

The center manifold $\{ u_- = \vec 0 \}$is invariant for the equation
\begin{equation}
\label{e:m:eq}
    \frac{d U}{d \tau} = F(U)
\end{equation}
and hence $F_- (\zeta, \, u_0, \, \vec 0) = \vec 0$ for every $\zeta$ and $u_0$. By regularity, 
$$
     F_-   (\zeta, \, u_0, \, u_-)  = G_s  (\zeta, \, u_0, \, u_-) u_-
$$
for a suitable matrix $G_s \in \mathbb M^{n_- \times n_-}$. Also, the uniformly stable manifold is invariant and hence proceeding as before we get that
$$
       F_0   (\zeta, \, u_0, \, u_-) = G_c (\zeta, \, u_0, \, u_-) u_0
$$  
for a suitable matrix $G_c \in \mathbb M^{n_0 \times n_0}$.  Finally, Hypothesis \ref{h:fast} implies that
$$
    f_1 (0, \,  \, u_0, \, u_-)   =0 
$$ 
and hence by regularity $ f_1 (\zeta, \,  \, u_0, \, u_-) = g_1 (\zeta,  \, u_0, \, u_-)  \zeta $.  Consider the decomposition
$$
    G_c (\zeta, \, u_0, \, u_-)  = G_{c } (\zeta, \, u_0, \, \vec 0) + \big[     G_c (\zeta, \, u_0, \, u_-) -   G_{c } (\zeta, \, u_0, \, \vec 0) \big]. 
$$
Thanks to Hypothesis \ref{h:center}, the subspace $\{ \zeta =0, \; u_- = \vec 0 \}$ is entirely constituted by equilibria and hence 
$$
    G_{c } (0, \, u_0, \, \vec 0) = \vec 0.
$$
By regularity, $G_c (\zeta, \, u_0, \, \vec 0) = G_{0 1} (\zeta, \, u_0) \zeta$ for a suitable matrix $G_{01} \in \mathbb M^{n_0 \times n_0}$ . Putting all the previous considerations together, we get that system \eqref{e:m:eq} can be written as 
\begin{equation}
\label{e:m:system1}
   \left\{
   \begin{array}{lll}
               \displaystyle{ d \zeta / d \tau =   g_1 (\zeta,  \, u_0, \, u_-)  \zeta  \phantom{\int} }\\
               \displaystyle{d u_0  / d \tau =  \Big\{    G_{0 1} (\zeta, \, u_0) \zeta + \big[  G_c (\zeta, \, u_0, \, u_-) -     G_c (\zeta, \, u_0, \, \vec 0) \big]  \Big\} u_0  } \\
               \displaystyle{ d u_- / d \tau =   G_s  (\zeta, \, u_0, \, u_-) u_- \phantom{\int}}\\ 
   \end{array}
   \right.
\end{equation}
Consider the decomposition
$$
     g_1 (\zeta,  \, u_0, \, u_-)   =  g_1 (\zeta,  \, u_0, \, \vec 0) + \big[  g_1 (\zeta,  \, u_0, \, u_-)   -  g_1 (\zeta,  \, u_0, \, \vec 0) \big]     
$$
By construction $G_s (0, \, \vec 0, \, \vec 0)$ admits only eigenvalues with strictly negative real part, thus ${G_s (\zeta, \, u_0, \, u_-) u_- = \vec 0}$ implies $u_- = \vec 0$. Thus, 
the set $\{ U: \; \zeta (U) = 0,  \;   F(U) = \vec 0 \}$ is the subspace $\{ \zeta =0, \; u_- = \vec 0 \}$. Thanks to Hypothesis \ref{h:slow}, we have 
$$
    g_1 (0,  \, u_0, \, \vec 0) = 0.
$$
By regularity, we thus have 
$$
      g_1 (\zeta,  \, u_0, \, \vec 0) = g_{11} (\zeta, \, u_0) \zeta 
      \qquad  
      \big[  g_1 (\zeta,  \, u_0, \, u_-)   -  g_1 (\zeta,  \, u_0, \, \vec 0) \big]   = G_{1 - } (\zeta,  \, u_0, \, u_-) u_-  
$$
for a suitable row vector $G_{1 - } (\zeta,  \, u_0, \, u_-) \in \mathbb R^{n_-}$.  Also, since the manifold $\{ u_0 = \vec 0, \; u_- = \vec 0 \}$ is entirely constituted by equilibria, then $g_{11}(\zz, \, \vec 0 ) = 0 $ for every $\zz$ and hence 
$$
      g_{11} (\zeta,  \, u_0) = G_{10} ( \zz, \, u_0 ) u_0 
$$ 
for a suitable vector $G_{10} \in \mathbb R^{n_0}$. 
In other words, \eqref{e:m:system1} reduces to 
\begin{equation}
\label{e:m:system2}
   \left\{
   \begin{array}{lll}
               \displaystyle{ d \zeta / d \tau =   \zeta^2 G_{10} (\zeta,  \, u_0 )   u_0 + \zeta   G_{1 - } (\zeta,  \, u_0, \, u_-) u_-  \phantom{\int} }\\
               \displaystyle{d u_0  / d \tau =  \Big\{    G_{0 1} (\zeta, \, u_0) \zeta + \big[  G_c (\zeta, \, u_0, \, u_-) -     G_c (\zeta, \, u_0, \, \vec 0) \big]  \Big\} u_0  } \\
               \displaystyle{ d u_- / d \tau =   G_s  (\zeta, \, u_0, \, u_-) u_- \phantom{\int}}\\ 
   \end{array}
   \right.
\end{equation}
\item \emph{Step 5:} we introduce a refined change of variables. 
Consider system 
\eqref{e:m:system1} restricted on the invariant subspace $\{ \zeta =0 \}$. One obtains 
\begin{equation}
\label{e:m:system3}
   \left\{
   \begin{array}{ll}
               \displaystyle{d u_0  / d \tau =  \Big[ G_c (0 , \, u_0, \, u_-) -     G_c (0, \, u_0, \, \vec 0)  \Big] u_0  } \\
               \displaystyle{ d u_- / d \tau =   G_s  (0, \, u_0, \, u_-) u_- \phantom{\int}}\\ 
   \end{array}
   \right.
\end{equation}
The subspace $\{ u_- = \vec 0 \}$ is entirely constituted by equilibria. Also, given a point $(u_0, \, u_-)$ belonging to a small enough neighbourhood of $\vec 0$, then 
the solution of \eqref{e:m:system3} starting at $(u_0, \, u_-)$ decays exponentially fast to a point in the subspace $\{ u_- = \vec 0 \}$. This is a consequence of the fact 
that $G_s (0, \, \vec 0, \, \vec 0)$ admits only eigenvalues with strictly negative real part. 

We can define a change of variables $\bar U =  \uu^4 (U)$ such that in the new variables $\bar U$ the following holds. For every $\bar u_0 (0) \in \mathbb R^{n_0}$ and for every $\bar u_- (0) \in \mathbb R^{n_-}$, the solution of  \eqref{e:m:system3} starting at the point $\big(\bar u_0(0), \, 
\bar u_-(0) \big)$ converges exponentially fast to the point $\big(\bar u_0(0), \, \vec 0 \big)$. In other words, the set $\{ \bar u_0 = \bar u_0 (0) \}$ is the stable manifold of system 
 \eqref{e:m:system3} around the equilibrium point $\big(\bar u_0(0), \, \vec 0 \big)$.  Let $\bar F (\bar U)$ be defined as in \eqref{e:m:barF}, with $\bar \uu = \uu^4$. Then 
$$
     F (\bar U) =
     \left(
     \begin{array}{ccc}
                   \bar  \zeta^2 \bar G_{10} ( \bar \zeta,  \,  \bar  u_0 )  \bar u_0 + \bar \zeta   \bar G_{1 - } (\bar \zeta,  \, \bar u_0, \, \bar u_-) \bar u_- \\
       \Big\{   \bar  G_{0 1} (\bar \zeta, \, \bar u_0) \bar \zeta + \big[  \bar G_c (\bar \zeta, \, \bar u_0, \, \bar u_-) -     \bar G_c (\bar \zeta, \, \bar u_0, \, \vec 0) \big]  \Big\} \bar u_0  \\
            G_s  (\bar \zeta, \, \bar u_0, \, \bar u_-) \bar u_-  \\
     \end{array}
     \right)
 $$
Because of the way we chose $ \uu_4$, when $\bar \zeta =0$ then $d \bar u_0 / d \tau = 0 $ and hence 
$$
      \big[  \bar G_c (0 , \, \bar u_0, \, \bar u_-) -     \bar G_c (0, \, \bar u_0, \, \vec 0) \big] \bar u_0 = \vec 0.
$$
By regularity, 
$$
     \big[  \bar G_c (\bar \zeta, \, \bar u_0, \, \bar u_-) -     \bar G_c (\bar \zeta, \, \bar u_0, \, \vec 0) \big]  = 
     \big[  \bar G_{0-} (\bar \zeta, \, \bar u_0, \, \bar u_-) -     \bar G_{0-} (\bar \zeta, \, \bar u_0, \, \vec 0) \big]  \zeta
$$
for a suitable function $G_{0 -} \in \mathbb M^{n_0 \times n_0}$.  
\item \emph{Step 6:} to conclude the proof of Proposition \ref{p:change} it is enough to define the local diffeomorphism $\uu$ as the composition of all the local diffeomorphisms defined in  the previous steps. 
\end{itemize}

\bibliography{biblio}
\end{document}